\documentclass[12pt]{article}
\usepackage[all]{xy}
\xyoption{textures}
\UseComputerModernTips

\newtheorem{theorem}{Theorem}[section]
\newtheorem{proposition}[theorem]{Proposition}
\newtheorem{lemma}[theorem]{Lemma}
\newtheorem{corollary}[theorem]{Corollary}
\newtheorem{definition}[theorem]{Definition}
\newtheorem{hypothesis}[theorem]{Hypothesis}
\newtheorem{remark}[theorem]{Remark}

\entrymodifiers={!!<0ex,.7ex>+}

\def\Om{\Omega}

\def\ac{{\cal A}}

\def\cc{{\cal C}}

\def\Z{{\bf Z}}

\def\la{\longrightarrow}
\def\be{\begin{equation}}
\def\ee{\end{equation}}

\def\lotimes{\stackrel{L}{\otimes}}

\def\oa{\stackrel{1}{+}}
\def\ob{\stackrel{2}{+}}
\def\oc{\stackrel{3}{+}}

\begin{document}

\title{Monoidal categories and multiextensions
\footnotetext{UMR  CNRS 7539,  Institut Galil\'ee,
Universit\'e Paris 13, F-93430 Villetaneuse.}
\footnotetext{ Email: breen@math.univ-paris13.fr}}
\author{Lawrence Breen}
\date{}
\maketitle

\begin{abstract}

We associate to a group-like monoidal groupoid $\cc$ a principal bundle $E$ satisfying
most of the axioms defining  a biextension. The obstruction to the existence of a genuine
biextension structure on $E$ is exhibited. When this obstruction vanishes,
the biextension $E$ is alternating, and a trivialization of $E$ induces a trivialization of
$\cc$. The analogous theory for monoidal  $n$-categories is also examined, as well as the
appropriate generalization of these constructions in a  sheaf-theoretic context.   In the
$n$-categorical situation, this produces a higher commutator calculus, in which some interesting 
generalizations of the notion of an alternating biextension occur. For $n= 2$, the corresponding
cocycles are constructed explicitly, by a partial symmetrization
process, from the cocycle describing the monoidal $2$-category.
\end{abstract}

\addtocounter{section}{-1}
\section{Introduction}

\setcounter{equation}{0}%

\hspace{.7cm}Let $A$ and $B$ be a pair of abelian groups. Central extensions 
\[ 0 \la A \la E \la B \la 0\]
of $B$ by $A$ are
classified up to equivalence by the $A$-valued cohomology group $H^2(B,A)$ (where $A$ is viewed as
a trivial $B$-module). To such an extension is associated the commutator map 
\[ \begin{array}{rcc}
\lambda: B \times B  & \la  & A \\
(b_1,b_2) & \mapsto &[s(b_1),s(b_2)]
\end{array}\]
determined by the choice of an arbitrary set-theoretic section $s$ of the projection from $E$ to
$B$. It is easily verified that this commutator map is independent of the
choice of the section $s$, and that it is a bilinear alternating map from $B \times B$ to $A$.
By construction, the map $\lambda$ measures the lack of commutativity of the group law of $E$. In
particular, the central extension $E$ is actually commutative whenever the map $\lambda$ vanishes,
so that it then determines an element of the group
$\mathrm{Ext}^1(B, A)$. 

\medskip
These facts, which are well-known, may be interpreted as follows in cohomological terms. Since the
group 
$B$ is abelian, its first integral homology group $H_1(B)$  is isomorphic to the group $B$
itself. Furthermore, the Pontrjagin product map
$H_1(B)
\times H_1(B) \la H_2(B)$ is bilinear, alternating, and  therefore induces a map
$\Lambda^2B \la H_2(B)$ which is an isomorphism. The previous discussion now follows directly by
considering the exact sequence
\be
\label{centr}
\mathrm{0} \la \mathrm{Ext}^1(B, A) \la H^2(B,A) \la \mathrm{Hom}(\Lambda^2 B, A) 
\ee
provided by the universal coefficient theorem. The functor $\Lambda^2$  here denotes the second
exterior power
$\Lambda^2_{\Z}$, applied to the group 
$B$ viewed as a $\Z$-module. Unless explicitly stated the corresponding higher exterior power 
functors
$\Lambda^j_{\Z}$ will in the sequel simply be denoted by $\Lambda^j$.

\medskip
Our aim in the present paper is to analyze in a similar manner some of  the higher cohomology groups
$H^n(B,A)$. These  have various geometric interpretations,  analogous to the
description of $H^2(B,A)$ in terms of central extensions, most of which are mentioned in
\cite{ml:hist}. The most general one of these interpretations of degree $n$ cohomology
groups provides a classification of
$n$-monoidal categories. 
 In the first case
of interest, that in which
$n=3$, this was first worked out  (for symmetric monoidal categories) in the barely accessible
\cite{sinh}, where $H^3(B,A)$ was interpreted as the group of equivalence classes of group-like 
monoidal groupoids
$\cc$, whose group
$\pi_0(\cc)$ of isomorphism classes of objects is isomorphic to $B$, and whose group $Aut_{\cc}(I)$
of self-arrows of the identity object $I$ of $\cc$ is isomorphic to the $B$-module $A$. 

\medskip
Our aproach
to the study of such monoidal category derives from the observation that there exists a natural
filtration, determined by powers of the augmentation ideal,  on the chains on a free abelian
simplicial resolution of 
$K(B,1)$. We  intend to study  this in some detail in \cite{lb:fh}, where we will examine the
effect of this filtration  on the integral homology of the abelian group
$B$. Let us merely  observe here that such a filtration on the chains of $B$ determines a
corresponding one on $A$-valued cochains, and therefore induces a filtration on the
 cohomology groups
$H^n(B,A)$. The {\it ith} associated graded piece of this filtration is  the group
$\mathrm{Ext}^i(L\Lambda^jB, A)$, where
$i+j = n$. Here  $L\Lambda^jB$ is the object in the derived category of abelian groups obtained by
applying the exterior power functor $\Lambda^j$ to a free abelian simplicial resolution of $B$
placed in degree zero, so that a more traditional notation for this  derived object would be 
$L\Lambda^j(B,0)$. 

\medskip
 Let us begin by considering  the case $n=3$. The filtration on $H^3(B,A)$  determines 
 three {\it a priori} non trivial terms in the associated graded group.
The first  of these is the group $\mathrm{Hom}(L\Lambda^3B, A)$. Since this group is isomorphic to  
$\mathrm{Hom}(\Lambda^3B, A)$,  a monoidal category $\cc$ of the type described above
determines a trilinear alternating map
$\varphi
\in
\mathrm{Hom}(
\Lambda^3 B, A)$. When this map $\varphi$ vanishes,  an element of the
group $\mathrm{Ext}^1(L\Lambda^2B,
A)$ can be associated to the category $\cc$. It was shown in  \cite{lb:alt} that this group
classifies, up to equivalence, the set of alternating biextension
$E$ of
$B
\times B
$ by
$A$. This may be understood in the present context by considering the commutator of $\cc$. This is
a principal
$A$-bundle
$E_{\cc}$ on
$B
\times B$, first introduced by P. Deligne in
\cite{pd:symbole}, whose fibre over an element $(x,y) \in B \times B$ is the set $E_{x,y}$ of
all arrows\footnote{We will
 generally denote by $XY$, rather than by the more customary $X
\otimes Y$, the product  of two objects $X$ and
$Y$ in a monoidal category 
$\cc$.}  $YX \la XY$ in
$\cc$ ($X$ and $Y$ being chosen representative objects in $\cc$ for the isomorphism classes $x$ and
$y$). Such a bundle may be endowed with a pair of partial composition laws determined by the
multiplication law in $\cc$, which are both associative, and  compatible with each other.
The obstruction to the commutativity of both of these partial group laws is described by  the
alternating map
$\varphi:
\Lambda^3B
\la A$ mentioned above. When 
$\varphi$ is trivial, the commutator $A$-bundle $E_{\cc}$ is therefore a genuine biextension of
$B\times B$ by
$A$, and in fact it is automatically an alternating one. Passing from the bundle $E_{\cc}$ to its
isomorphism class, we  obtain in this manner the sought-after element of the group
$\mathrm{Ext}^1(L\Lambda^2B, A)$.  When this element vanishes, the category $\cc$
determines an element in the last component of the graded group associated to $H^3(B,A)$, in other
words in the group
$\mathrm{Ext}^2(B,A)$. In geometric terms, this may be interpreted as the assertion that a
trivialization of
$E_{\cc}$ as an alternating biextension determines on
$\cc$ a strictly symmetric monoidal structure. By
\cite{pd:sga4}, we know that such strictly symmetric group-like monoidal groupoids are
indeed classified by the sought-after group
$\mathrm{Ext}^2(B, A)$. However, this group of extensions always  vanishes in the
category of abelian groups,  so that it does not provide a genuine  invariant attached to
$\cc$. In geometric term, this is reflected in the assertion that such a  strictly symmetric
group-like monoidal groupoid is always  equivalent to the trivial symmetric monoidal
category associated to the pair of groups
$B$ and $A$. 

\medskip
It is instructive to carry out the previous discussion purely in terms of a given $A$-valued
3-cocycle $f(x,y,z)$ on $B$. The alternating map $\varphi$
which one then encounters is a very familiar one, being simply the map obtained by evaluating $f$ on
the decomposable elements of
$H_3(B)$. In order to interpret the commutator biextension $E_{\cc}$ directly in terms of
$f(x,y,z)$, we have found it necessary to insert in our text a description  of alternating
biextensions in purely cocyclic terms. We believe that such a description, which was not
carried out in
\cite{lb:alt},  can be of independent interest. It turns out  that the pair of cocycles
($g(x,y;z), h(x;y,z)$) which  describe the commutator biextension $E$  of
$\cc$ are obtained from
the given 3-cocycle $f(x,y,z)$ by a  partial symmetrization process which  already occurs (without
the assumption that $B$ is  abelian)  in a computation by R. Dijkgraaf and E. Witten of the
2-cocycle associated by Chern-Simons theory  to a  given $3$-cocycle  \cite{di-wi} \S 6.6.

\medskip
The rest of this text is devoted to various generalizations of the previous discussion. The
first of these  extends the theory from the study of monoidal categories to that of monoidal
stacks. This level of generalization is analogous to that which occurs when one passes from the
classification of central extensions of abelian groups to that of topological abelian groups
\cite{gs:topgroups} or of algebraic groups (\cite{gacc} chapter VII). The  choice of  objects or
arrows in $\cc$ required for a cocyclic description of the monoidal stack $\cc$ can in general only
be made locally. The cohomological obstruction to a global choice of  objects is determined, as
explained in
\cite{lb:2-gerbes}, by the underlying gerbe of  $\cc$. The obstruction
to a corresponding global choice of arrows is reflected in the fact that the commutator biextension
$E_{\cc}$ of  $\cc$ no longer has, as in the category case, a global section above its base $B
\times B$.
$E_{\cc}$ is  now a genuine biextension of $B \times B$ by $A$ in the sense of \cite{ag:sga7},
rather than  one which may be described, as in the category case, by a pair of cocycles $(g,h)$.

\medskip
Our next generalization consists in passing from the cohomology group 
$H^3(B,A)$ to the group $H^4(B,A)$. The latter  classifies the monoidal 2-groupoids $\cc$
which satisfy the conditions $\pi_0(\cc)= B,\pi_1(\cc)= 0$ and $
\pi_2(\cc)= A$. The natural action of $B$ on
$A$ is once more assumed to be trivial. A geometrical discussion of the associated graded pieces
for the filration on $H^4(B,A)$  requires a geometrical  understanding of the corresponding groups
$\mathrm{Ext}^i(L\Lambda^jB , A)$ for 
$i + j = 4$ . We interpret these  groups as the groups of equivalence classes of certain geometric
objects which we call the 
$(i,j)$-extensions of
$B$ by
$A$. When $i=1$, these are simply, for an arbitrary $j$,  the $j$-fold extensions
of
$B$ by
$A$ introduced by A. Grothendieck in \cite{ag:sga7}. We  will therefore use this concept here for
$j=3$, and  we will call such objects triextensions of $B$ by $A$. The next term in the filtration
requires that we understand the notion of a
$(2,2)$-extension of $B$ by $A$. This is an interesting new concept, consisting in
a category (or more generally a stack) 
$\mathcal{E}$ for which $\pi_0(\mathcal{E}) = B \times B$ and $\pi_1(\mathcal{E}) = A$, and which
is endowed with a pair of coherently associative and appropriately compatible partial group
laws, which define on the restrictions of
$\mathcal{E}$ to all subsets $x\times B$ and
$B\times y$ the structure of a group-like symmetric monoidal category.

\medskip

While these definitions of  a triextension and of a (2,2)-extension  present no
great difficulty, there remains the question of imposing on each of these objects an alternating
structure. In order to achieve this, we make use of Koszul complex techniques, and interpret
the requisite groups
$\mathrm{Ext}^i(L\Lambda^jB , A)$ in geometric terms. Once the appropriate definitions have
been obtained,  we can describe the geometric objects which our higher commutator
calculus associates to a given monoidal 2-category. Part of this discussion is carried out in
cocyclic terms, an efficient substitute in the present context for  pasting diagrams
in  2-categories.  A pleasant feature of this discussion is the occurence of a systematic partial
symmetrization process, analogous to  the one mentioned above in 3-cocycle situation, and which
points quite clearly to a general statement for the corresponding filtration on the cohomology
groups of
 arbitrary degree.  Certain of these symmetrized higher cocycles occur, for a non-abelian group $B$,
as the images of higher transgressions in recent work of  J.-L. Brylinski and  D. A. McLaughlin
\cite{br-ml}. 

\medskip

We have assumed throughout this text that the group $B$ was abelian, but the constructions carried
out here remain for the most part valid without that hypothesis. Indeed, in  Deligne's original
construction \cite{pd:symbole} of the commutator $E_{\cc}$ of a  monoidal category $\cc$, no such
commutativity assumption on the group law of $B$ was made, nor was it required in the
previously mentioned texts \cite{di-wi} and \cite{br-ml}. Without such a commutativity hypothesis,
the torsor
$E_{\cc}$ is only defined above that part of $B \times B$  which consists of pairs  of
commuting elements $x,y$ of
$B$. The partial group laws which are introduced here  no longer yield in that case a  biextension,
but a weaker structure which deserves to be formalized. While we have not carried out this
formalization here in order not to overburden this text, we intend to return to this question in the
future. Let us simply observe for the present that those central extensions whose associated
commutator maps are the most interesting are central extensions for which the quotient group $B$ is
not abelian. It is therefore to be expected that the same will be true for the higher
constructions which we examine here.
  
\medskip

While we have  emphasized in this introduction the cohomological interpretation of our
constructions, in terms of the derived functors of the exterior algebra functor, this will not be
the case in the sequel.
 Indeed, the emphasis will henceforth be on the determination of the new higher alternating
structures, rather than on the quest for an interpretation in geometric terms of the 
universal coefficient theorem. This text is therefore independent of the forthcoming
\cite{lb:fh}. Both approaches are, however, fully compatible, and shed light upon each other. In
the present context, this is illustrated  in 
\cite{lb:alt}, where   alternating biextensions are analyzed {\it via} the universal coefficient
theorem.

\section{Cohomology and categories}
\label{sub:uct}

\setcounter{equation}{0}%

\hspace{.7cm}The most general  interpretation of  the cohomology group $H^3(B,A)$ is the one
 due to A. Grothendieck. It expresses degree three cohomology classes  in terms of monoidal
categories (see
\cite{sinh}, \cite{ca-ce} \S 2.1
 and also,
 in a sheaf-theoretic context in which the 3-cocycles do not appear explicitly,
  \cite{pd:sga4}). 
 We begin by recalling this interpretation of  $H^3(B,A)$, and refer to 
\cite{kbrown} IV \S 5,  and to \cite{ml:hist} and references therein for related descriptions
of this cohomology group.
 Observe first of all that if $({\cal C},\otimes , a)$ is a monoidal group-like 
groupoid\footnote{also referred to as a {\it gr}-category \cite{sinh} or a categorical
group \cite{js:braided}.} with
unit object $I$, then the monoidal structure on 
${\cal C}$ determines, for each object $X \in
\cc$,  a right multiplication isomorphism
\begin{equation}
\label{autb}
\xymatrix{
{A = Aut(I)} \ar[r]^<<<<<{\otimes X}& {Aut(X)}}
\end{equation}
through which any group of automorphisms in $\cc$ will
henceforth be identified with $A$. The image in $Aut(X)$ by the left
multiplication isomorphism
\be
\label{autb1}
A = Aut(I) \stackrel{X \otimes}
\longrightarrow Aut(X)
\ee
of an element $a \in A$ may be identified by (\ref{autb}) with an
element $^{X}a \in A$ which actually only depends on the isomorphism
class $x$ of $X$ in the group $B$ of isomorphism classes of objects of $\cc$, and which will therefore
be denoted
${}^{x}a$. It is readily verified that this action of $B$ on $A$ endows $A$ with a
$B$-module structure. 

\medskip
We now choose, for each  $x \in B$, an object $X_x \in \cc$  whose isomorphism class is
$x$. For each pair of elements $x, y \in B$, the objects $X_xX_y$ and $X_{xy}$  both live in the
component of $\cc$ described by the element $xy \in B$, so that there exist arrows between
them. Choose such an arrow
\be
\label{cxy}
c_{x,y}: X_x X_y \la X_{xy}
\ee
 for each $x,y \in B$.
For every $x,y,z \in B$, 
the associativity isomorphism 
\[ a_{x,y,z}: X_x(X_yX_z) \la (X_xX_y)X_z
\]
 determines an
element $f(x,y,z) \in \mathrm{Aut} (X_{xyz}) = A$ such that the  diagram

\be
\label{diag:f}
\xymatrix{
 {X_x(X_yX_z)} \ar[r]^{a_{x,y,z}} \ar[d]_{{X_x}c_{y,z}}& {(X_xX_y)X_z} \ar[d]^{c_{x,y}X_z}\\
{X_x X_{yz}} \ar[d]_{c_{x,yz}} & {X_{xy} X_z}  \ar[d]^{c_{xy,z}}\\
{X_{xyz}} \ar[r]_{f(x,y,z)} &
{X_{xyz}}}
\ee
commutes.
The pentagon axiom in $\cc$ then implies that $f(x,y,z)$ is a 3-cocycle. We may even assume, by
choosing the objects $X_x$ and the arrows (\ref{cxy}) carefully, that the 3-cocycle $f$ is
normalized (as will be all those occuring from now on, unless explicitly stated). Other
choices for these  objects and arrows  of $\cc$ yield a cohomologous
3-cocycle, so that the class of $\cc$ in $H^3(B, A)$, for the $B$-module structure on $A$
determined by (\ref{autb})-(\ref{autb1}), is well-defined.

\begin{remark}
\label{rem13} {\rm {\it i}) The previous construction may be interpreted as follows in
topological terms.
  The nerve $\mathcal{G} = N\cc$ of $\cc$ is a 2-stage
Postnikov system with homotopy groups
$\pi_{0}(\mathcal{G}) = B$  and $ \pi_{1}(\mathcal{G} , I) = A$. The monoidal structure on $\cc$
determines an
$A_\infty$ 
$H$-space structure on $\mathcal{G}$, so that  $\mathcal{G}$ deloops to a connected space
$\mathcal{X}=B\mathcal{G}$ whose homotopy groups
 $B$ and $A$ live respectively in degrees one and two. The $k$-invariant $k \in H^3(B,A)$ of 
the two-stage system $\mathcal{X}$ is the sought-for cohomology class associated  to the monoidal
category
$\cc$. Conversely, one can start from such 2-stage Posnikov system $\mathcal{X}$. 
The space $\mathcal{Y} = \Omega \mathcal{X}$ of loops on $\mathcal{X}$ is essentially the nerve of a
groupoid
$\cc$, and the
$H$-space structure  on $\mathcal{Y}$ corresponds to the monoidal structure on $\cc$.

\medskip
\indent {\it ii)} When the group law in $\cc$ is strict, the monoidal category $\cc$
 may be represented by a crossed module $N \la E$ with $E$ ({\it resp.} $N$)  the group of
objects ({\it resp.} the group of arrows sourced at the identity object) of $\cc$. A comparison
between the terms in the  formula \cite{kbrown} IV (5.7) and  the arrows in   diagram (\ref{diag:f})
implies that the  standard  method \cite{kbrown} IV \S 5 for associating
a 3-cocycle to a crossed module is consistent with the one  given above.

\medskip
\indent {\it iii)} Suppose that the monoidal category $\cc$ is endowed with a  commutativity
isomorphisms
$s_{x,y}: X_y X_x
\la X_x X_y
$. This determines, via the commutative diagram
\[
\xymatrix{
{X_y X_x} \ar[r]^{s_{x,y}} \ar[d]_{c_{y,x}} & {X_x X_y}\ar[d]^{c_{x,y}}
\\
{X_{xy}}
\ar[r]_{g(x,y)}&{X_{xy}}}
\]
a map 
\[
g: B \times B \la A
\]
for which the braiding axioms \cite{js:braided}
imply that the cocycle conditions
\begin{eqnarray}
\label{eq:braid}
f(x,y,z) - f(x,z,y) + f(z,x,y) & = & g(x+y,z) - g(x,z) - g(y,z) \\
- f(x,y,z) + f(y,x,z) - f(y,z,x) & = & g(x,y+z) - g(x,y) - g(x,z) \nonumber
\end{eqnarray}
are satisfied. Alternate choices yield a well-defined class in $H^4(K(B,2), A)$
\cite{js:braided} proposition 3.1,  \cite{lb:2-gerbes} \S 7.8. The braiding axioms allow  a
double delooping of the nerve $\mathcal{G}$ of $\cc$ to a two-stage space $\mathcal{Y}$, whose
$k$-invariant is this cohomology class. When
$\cc$ is symmetric monoidal, the
additional condition
\[
g(x,y) = g(y,x)
\]
is satisfied. The two conditions (\ref{eq:braid}) then coalesce and a class in the stable
cohomology group $H^5(K(B,3), A)$, corresponding to the $k$-invariant of a further delooping of
$\mathcal{G}$, is defined. It follows that
$\mathcal{G}$ is in that case an infinite loop space. Finally, when the stronger condition 
\be
\label{eq:ps}
g(x,x) = 0 
\ee
is satisfied, the monoidal category $\cc$ is strict Picard \footnote{in other words
a strictly symmetric 
 group-like monoidal groupoid.}. Prolonging by one step the canonical
resolution \cite{ag:sga7} VII (3.5.1)  of the abelian group $B$, it is apparent that the pair
$(f,g)$ now describes a class in the  group $\mathrm{Ext}^2(B, A)$. The nerve of $\mathcal{G}$ is
now equivalent to a simplicial abelian group, and the appropriate trunctation of its associated
Moore complex determines the class in question. We observed earlier that such a group 
$\mathrm{Ext}^2$  is always trivial in the category of abelian
groups. All such strict
Picard categories are therefore  equivalent to trivial ones. The cocycle
$(f,g)$ which describes such a  strict Picard  category $\cc$ is therefore a coboundary, so that
there exists a map $h: B^2
\la A$ satisfying the following conditions:
\begin{eqnarray}
\label{eq:pstriv}
  f(x,y,z) & =  &  h(y,z) - h(x+y, z) + h(x,y+z) - h(x,y)\\
 g(x,y) & = & h(x,y) - h(y,x)
\end{eqnarray}
In particular the underlying monoidal category of $\cc$ is trivial, as  reflected by the fact
that $h$ is a group cohomology coboundary for $f$.
 }
\end{remark}

We  now investigate the conditions under which the given group law in
the monoidal category $\cc$ satisfies some form of commutativity. One possible approach would
consist in expressing in geometric terms  the  obstructions to the surjectivity
of the successive suspension maps
\[ H^5(K(B,3), A) \la H^4(K(B,2), A) \la H^3(B, A) \]
which control the level of commutativity of the group law on $\cc$. We will instead examine here,
without passing through these intermediate steps, the conditions under which the monoidal
category    $\cc$ can be endowed with a fully symmetric monoidal structure.  Let us begin by making 
 the following additional assumption.
\begin{hypothesis}
\label{hyp1}
The group $B$ is abelian, and the
$B$-module structure on $A$  is trivial.
\end{hypothesis}

Note that this is a very weak commutativity condition.
Indeed, the requirement that  there exists, for each pair of objects $X,Y$ in $\cc$,
 an isomorphism
between $YX$ and $XY$
implies that the group $B$
of isomorphism classes of objects of $\cc$ is abelian. If we also ask
that this family of isomorphisms be natural in the objects $X$ and $Y$,
 and compatible in the obvious sense with the
identity object, then  the $B$-module structure on $A$
is trivial, the triviality being
expressed by the commutativity of the following diagram
\[
\xymatrix{& {X} 
\ar[dl] \ar[dr] & \\{XI} \ar[d]_{Xf} \ar[rr] && {IX } \ar[d]_{fX}\\{XI} 
\ar[rr]  \ar[dr] &&{IX} \ar[dl] \\& {X} &}
\]
In particular, hypothesis \ref{hyp1} is automatically satisfied whenever the
category
$\cc$ is braided. 

\medskip

Hypothesis \ref{hyp1} allows us to apply the universal coefficient theorem to the computation of 
$H^3(B,A)$, and therefore to obtain an
analog for $H^3$ of the exact sequence (\ref{centr}). Let us  begin with the naive approach to
this question.   The terms of the universal coefficient exact sequence
\be
\label{uct3}
0 \la \mathrm{Ext}^1(H_2(B), A) \la H^3(B,A) \la \mathrm{Hom}(H_3(B), A) \la 0
\ee
can be made explicit since  the appropriate homology groups are known (see \cite{lb:fh}). To
the given class in  $H^3(B,A)$  of  a monoidal group-like groupoid  $\cc$ satisfying hypothesis
\ref{hyp1} is associated a trilinear alternating map 
$\varphi_\cc : \Lambda^3B \la A$. If
$\varphi$  vanishes, then a second map $\psi_{\cc}: \Om_2 B \la A $ may be associated to $\cc$.
The source $\Om_2B$ of this arrow is a  group $\Om B$  first defined by Eilenberg-Mac Lane
\cite{eml:hpin}, and which can be interpreted as the first (non-additive) left derived functor
$L_1\Lambda^2(B,0)$ of the abelian group $B$ set in degree zero. A specific presentation
presentation  of
$\Om_2 B$ expresses 
the map $\psi_{\cc}$ in terms of a family of quadratic maps $\psi_n : {}_nB \la A $ for
varying positive integers $n$, related to each other in an appropriate manner. Finally, if
$\psi_{\cc}$ is trivial, the class of $\cc$ is  determined by an element $\chi_{\cc}$ in the
extension group
$\mathrm{Ext}^1(\Lambda^2B, A)$.

\medskip
While this is a complete discussion, there remains the question of interpreting it in geometric
terms. The fact that the functor $\Om_2 B$ is the first derived functor of the exterior
algebra functor $\Lambda^2B$ suggests that the real object of  interest, encompassing both
$\psi_{\cc}$ and  $\chi_{\cc}$, lives in the group
$\mathrm{Ext}^1(L\Lambda^2B, A)$. Indeed it is shown in \cite{lb:alt} remark 3.6 that an element
in this group determines by d\'evissage appropriate elements $\psi_{\cc}$ and  $\chi_{\cc}$. The
next three sections will provide a construction of the object
$E_{\cc}$  whose class determines the sought-after element in the group
$\mathrm{Ext}^1(L\Lambda^2B, A)$.

\section{Alternating biextensions}

\setcounter{equation}{0}%

\hspace{.7cm}The group $\mathrm{Ext}^1(L\Lambda^2B, A)$  was given a geometrical interpretation
in \cite{lb:alt}, as the group of equivalence classes of alternating biextensions of $B$ by
$A$. Let us begin by  reviewing the definition of  an alternating biextension.
Let  $E$ be an $A$-torsor above $B \times B$. Its fiber above a point $(x,y) \in B \times B$ will
be denoted  $E_{x,y}$.  Recall first of all that an (ordinary) biextension  of $B \times B$ by $A$
is such  an
$A$-torsor
$E$ above $B \times B$, endowed with a pair of
partial composition laws whose restrictions to the appropriate fibers may be depicted as morphisms
of $A$-torsors 
\begin{equation}
\label{eq: mult1}
\stackrel{1}{+} :E_{x,y} \wedge E_{x',y}
\longrightarrow E_{xx',y} \label{eq:+1} 
\end{equation}
\be
\label{eq: mult2}
\stackrel{2}{+} :E_{x,y} \wedge E_{x,y'}
\longrightarrow E_{x,yy'} \label{eq:+2}   
\ee 
where $\wedge = \wedge^{A}$ denotes the contracted product of the corresponding $A$-torsors. These
two composition laws are required to be associative, commutative and compatible with each other
\cite{dm:biextensions},
\cite{ag:sga7} expos\'e VII \S 2.
A torsor endowed with a pair of partial multiplication laws which are merely associative and
compatible will be called a weak biextension\footnote{In that case, a better notation for 
the two  partial group laws would be
$\stackrel{1}{\times}$ and $\stackrel{2}{\times}$.}.

\medskip
We review
for the reader's convenience the manner in which the structure on a biextension $E$
 whose underlying torsor is trivialized may be described in terms of  cocycles.
The triviality hypothesis asserts that the underlying torsor $E$ may simply be defined by 
$E =  A \times B \times B $.
The first and second partial group laws are respectively determined by  maps
$g(b_1, b_2 ; b')$ and $h(b; b'_1, b'_2)$ from $B^3$ to $A$ such that 
\be
\begin{array}{ccccc}
\label{ghdef}
(a, b_1, b') & \stackrel{1}{+} & (a, b_2, b') & = & (a + g(b_1, b_2; b'), \: b_1 + b_2,\:
b') \nonumber
\\ (a, b, b'_1) & \stackrel{2}{+} & (a, b, b'_2) & = & (a + h(b; b'_1, b'_2) , \: b, \: b'_1
+ b'_2)
\end{array}
\ee
on $E$. The associativity conditions for these laws 
 translate to the following
 cocycle conditions\footnote{ The group law of $B$ will henceforth be written additively,
and that of
$A$ multiplicatively.} on $g$ and $h$ \cite{dm:biextensions}.
\begin{eqnarray}
\label{bicocyc}
\frac{g(b_2, b_3 \, ; b') \:  g(b_1, b_2 + b_3 \, ; b')}{ g(b_1 + b_2, b_3 \, ; b') \: g(b_1, b_2
\, ; b')}  & = & 1 \\
&& \nonumber\\
\frac{h(b \, ; b'_2, b'_3) \:  h(b \, ; b'_1, b'_2 + b'_3)}{  h(b \, ; b'_1 + b'_2, b'_3)  \: h(b
\,; b'_1, b'_2)} & = &1 \nonumber
\end{eqnarray}
in other words to the standard 2-cocycle condition for the maps
$g( - \: , \: - \,  ; \, b')$ and $h(b \,; \: - \, , \,  -)$  from $B^2$ to $A$, for all fixed 
$b, b' \in B$. Similarly, the commutativity conditions, when they are satisfied, translate to the
standard symmetry cocycles with the last ({\it resp.} the first) variable fixed
\begin{eqnarray*}
\label{bicom}
g(b_1, b_2 \, ; b') & = & g(b_2, b_1 \, ; b') \\
h(b \, ; b'_1, b'_2) & = & h(b \, ; b'_2, b'_1) \nonumber
\end{eqnarray*}
Finally, the compatibility condition now becomes the rule
\be
\label{bicompat}
\frac{h(b_1 + b_2\, ; b'_1, b'_2)}{  h(b_1 \, ; b'_1, b'_2) \:  h(b_2 \, ; b'_1, b'_2)} = 
\frac{g(b_1, b_2 \, ; b'_1 + b'_2)}{g(b_1, b_2 \, ; b'_1) \: g(b_1, b_2 \,  ; b'_2)}
\ee 
 A  cocycle pair
$(g,h)$ is cohomologous to zero, and therefore defines a trivial
 biextension
structure, whenever there exists a map $k: B \times B \la A$ such that
\begin{eqnarray}
\label{bitriv}
 g(b_1, b_2 \, ; b') & =  &\frac{ k(b_1 + b_2 , b')}
{ k(b_1 ,  b') \: k(b_2 ,  b')} \nonumber \\
&&  \\
 h(b \, ; b'_1, b'_2)  & = & \frac{ k(b, b'_1 + b'_2) }
{ k(b, b'_1) \:  k(b, b'_2)} \nonumber 
\end{eqnarray}

 \medskip
Let us now pass from ordinary  to alternating biextensions. Their description in \cite{lb:alt}
was modelled on the exact triangle derived from the Koszul sequence
\be
\label{kosl2}
\mathrm{0} \la \Gamma_2B \la B \otimes B \la \Lambda^2B \la \mathrm{0}
\ee
One begins by considering an (ordinary) biextension $E$ of $B \times B$ by $A$. The restriction
$\Delta E$ of  $E$ to the diagonal in $B \times B$ is an $A$-torsor on $B$
 with the following additional properties.
\begin{enumerate}
\item The torsor  $\Delta E$ is symmetric, in other words
there exist,
 for
 varying $x \in B$, 
a family of symmetry isomorphisms
\be
\label{symiso}
 \sigma _x: \Delta E_{-x} \la \Delta E_x 
\ee
\item  $\Delta E$ is  endowed with a cube structure which is 
 compatible, in a strong sense, with the symmetry isomorphism.
\end{enumerate}
 The precise sense in which the cube structure on the $A$-torsor
 $L = \Delta E$ on $B$ is
compatible
with the symmetry isomorphism $\sigma:i^* L \la L$ (\ref{symiso}) is best explained as follows.
It is essentially shown in 
 \cite{lb:cube} \S 5,  though not made explicit there, that for any $A$-torsor $L$ on $B$ endowed
with a cube structure  the
$A$-torsor
$L \wedge i^*L^{-1}$   on $B$ is canonically endowed with a composition law, which makes it into a
(commutative) extension of $B$ by $A$. The  compatibility between cube structure and symmetry
on $L$ may be expressed as the requirement that the symmetry isomorphism $\sigma$, viewed as a
section  of
$L
\wedge i^*L^{-1}$ on $B$,  splits it as a group  extension. When this condition is satisfied, one
says that the
$A$-torsor $L$ is endowed with a 
$\Sigma$-structure. An alternating structure on a biextension 
$E$ is then defined as follows. 
\begin{definition}
\label{defalt}
An alternating biextension of $B \times B$ by $A$ is a biextension $E$ of $B \times B$ by $A$,
 together with a trivialization $t: B \la E$ of the restriction $\Delta E$ of $E$ to the diagonal
compatible  with the $\Sigma$-structure of $\Delta
E$.
\end{definition}

  When the underlying torsor of $E$ has a global section, this definition of an alternating
biextension can be made explicit
 in terms of the pair of cocycles $(g,h)$ attached to $E$.
 A trivialization $t$ of $\Delta E$ compatible with the symmetry isomorphism (\ref{symiso}) is
expressed by  a map
$u: B
\la A$
 (for which we may assume that $u(0) = 1$)
such that
\be
\label{symcoc1}
 \frac{u(-b)}{u(b)} = \frac{g(b,-b \, ;b)}{h(-b \, ; b, -b)} 
\ee
for all $b \in B$. We now introduce a map $\lambda: B^2 \la A$
which may, in view of (\ref{bicompat}), be defined by either of the two following equations 
\begin{eqnarray}
\label{lambdabiex}
 \lambda (b_1, b_2) & =
 &  g(b_1, b_2 \, ; b_1 + b_2) \:  h(b_1\, ;  b_1, b_2) \: h(b_2 \, ; b_1, b_2) \nonumber \\
 &= & h(b_1 + b_2 \, ; b_1, b_2) \:  g(b_1, b_2 \, ; b_1) \: g(b_1, b_2 \, ; b_2)
\end{eqnarray}
The requisite compatibility between the trivialization 
$u$ of $\Delta E$ and the cube structure on $\Delta E$ may now be expressed as the
condition
\be
\label{symcoc2}
\Theta (u)(b_1, b_2, b_3) = \frac{\lambda(b_1 + b_2, b_3)}
{\lambda (b_1, b_3) \lambda (b_2, b_3)} \:
g(b_1, b_2 \, ; b_3) \:  h(b_3 \, ; b_1, b_2)
\ee
where  $\Theta (u)$ is   the second difference of the
map
$u$, defined by
\[
\Theta (u)(b_1, b_2, b_3) = \frac{u(b_1 + b_2 + b_3) \, u(b_1) \, u(b_2) \, u(b_3)}
{u(b_1 + b_2) \, u(b_1 + b_3) \,  u(b_2 + b_3)}
\]
Finally,  an alternating biextension $(g,h,u)$ is trivialized by a map $k: B \times B
\la A$ satisfying the trivialization conditions (\ref{bitriv}), together with the additional
condition
\begin{eqnarray}
\label{ktriv}
k(x,x) & = & u(x)
\end{eqnarray}

\medskip
A somewhat more intuitive description of an alternating biextension is
obtained by introducing first  the simpler concept of an anti-symmetric biextension. Consider the
functor
$\mathrm{As}^2B$ of anti-symmetric tensors on $B$, which fits into the following commutative
diagram whose
 horizontal lines are exact when $B$ is free.
\[
\xymatrix{
{0} \ar[r]  & {\mathrm{Sym}^2B} \ar[r] \ar[d] & {\otimes^2B} \ar[r] \ar[d] & {\mathrm{As}^2B} \ar[r]
\ar[d] & {0}
\\{0} \ar[r] 
&{\Gamma_2B}\ar[r] 
& {\otimes^2B}\ar[r]
&{\Lambda^2B}\ar[r]
&{0}}
\]
By the snake lemma, this determines for every free abelian group $B$ a short exact sequence
\be
\label{as2}
0 \la B/2B \la \mathrm{As}^2B \la \Lambda^2B \la 0
\ee
 Anti-symmetric biextensions are to be thought of as those
biextensions
$E$ of
$B
\times B$ by
$A$ which are classified up to isomorphism by the group
$\mathrm{Ext}^1(L\mathrm{As}^2B, A)$. Denoting by $s$ the map from $B^2$ to itself which permutes
the factors, this means that they are the biextensions $E$ for which we are given  a trivialization
$\pi$ of the induced
 biextension $F = E \wedge s^*E$ compatibly with the natural symmetric biextension structure on
$F$. Such a trivialization $\pi$  may  also be described by the induced 
 biextension isomorphism
\be
\label{pialt}
\pi_{x,y}: E_{x,y}^{-1} \la  E_{y,x} 
\ee
between $E^{-1}$ and the pullback $s^*E$ of $E$.
The symmetry condition on $\pi$ then becomes the requirement that for each $(x,y) \in B \times B$, 
the map
${}^t\pi_{x,y}:E_{y,x}^{-1}
\la E_{x,y}$ induced by
$\pi$ coincides with   $\pi_{y,x}$. It is readily verified  that any alternating biextension is
anti-symmetric (\cite{lb:alt} proposition 1.4). The distinguished triangle associated to
(\ref{as2}) gives us  a new description of an alternating biextension $E$ in terms of the underlying
anti-symmetric one. Observe first of all that for any anti-symmetric biextension
$E$ of $B \times B$ by $A$, the pullback 
$\Delta E$ of
$E$ along the diagonal is actually a commutative extension\footnote{The
 definition of the group law on $\Delta E$ given in (\ref{ustarv}) below for a particular
biextension
$E$ is valid in the general case.} of $B$ by $A$. Furthermore, the pullback of this extension by the 
``multiplication by
$2$'' map 
\be
\label{mult2}
2_B: B \la B
\ee
 on $B$ is canonically split as an extension. This is equivalent to the assertion that the square
$\Delta E^2$ of the extension $\Delta E$ (under Baer addition) is split.  An alternating
biextension may  now be described in the following manner.

\begin{proposition}
\label{altantisym}
An anti-symmetric biextension $E$ is alternating if and only if its restriction
$\Delta E$ to the diagonal is  split as an extension, by a  splitting which is compatible with the
 splitting of $\Delta E^2$ determined by the anti-symmetry structure on $E$.
\end{proposition}

Here is the cocyclic translation of this new description of an alternating biextension, when the
underlying torsor of the biextension $E$ is trivial. The biextension structure on such a
trivial $A$-torsor is described, as before, by  a pair of maps  
$g(b_1, b_2 \, ;\:  b')$ and $h(b\, ;\:  b'_1, b'_2)$ (\ref{ghdef}). An anti-symmetry structure on
$E$ is determined by a map $\varphi : B^2 \la A$  which trivializes the induced biextension $F$, in
other words a map $\varphi$ such that the equations
\begin{eqnarray}
\label{phitriv}
\frac{\varphi (b_1 + b_2, b')}{\varphi(b_1, b') \:  \varphi(b_2, b')}& = & g(b_1, b_2 ; b') \: 
h(b';  b_1, b_2) \nonumber \\
\frac{\varphi (b, \, b'_1 + b'_2)}{\varphi (b, b'_1) \: \varphi(b, b'_2)}& = & h(b; b'_1, b'_2) \:
g(b'_1, b'_2; b)
\end{eqnarray}
are satisfied. Since the trivialization of $F$ defined by $\varphi$ must be compatible with the
symmetry structure on $F$, the maps $\varphi $ must satisfy the additional condition
\be 
\label{phisym}
\varphi (b,b') = \varphi (b', b)
\ee
for all $ b, b' \in B$. The following assertion is proved by a rather elaborate cocycle
computation, which we omit.
\begin{lemma}
The map
$c: B^2 \la A$  defined by 
\be
\label{defc}
c(b, b') = \lambda (b, b') \: \varphi (b,b')
\ee
(where $\lambda (b,b')$ is given by (\ref{lambdabiex})) is an $A$-valued 2-cocycle on $B$.
\end{lemma}

 The
commutativity condition for the partial group laws $\stackrel{1}{+}$ and $\stackrel{2}{+}$,
together with equation (\ref{phisym}), imply that the 2-cocycle $c$ is symmetric, so that the
triple 
$(g, h,
\varphi)$  determines a commutative extension of $B$ by $A$, which in  fact  is one
previously obtained by  restricting
 the anti-symmetric biextension $E$ above the diagonal in $B \times B$. The equations
(\ref{bicocyc}) - (\ref{bicompat}) imply that 
\[ 
 \lambda (b, b')^2 = \frac{\varphi (b+b', b+b')}{\varphi (b,b) \: \varphi (b,b') \:
\varphi (b',b) \: \varphi (b',b')}
\]
so that the equation
\[
c(b,b')^2 = \frac {\varphi(b+b', b+b')}{\varphi(b,b) \: \varphi (b', b')}
\]
is satisfied. This shows that the  1-cochain  $\psi(b)$ defined by 
\[ \psi (b) = \varphi (b, b)
\]
 trivializes the 2-cocycle $c(b,b')^2$ which describes $\Delta E^2$. Taking into account
the significance of a trivialization of $E$ compatible with all this structure, we may now summarize
the previous discussion in the following way.
\begin{proposition}
\label{quadrup}  An alternating biextension of $B \times B$ by $A$, with trivial underlying
torsor, is determined by a quadruple 
\be
\begin{array}{cccc}
{\xymatrix{
{B^3} \ar[r]^g &
{A}}}&
{\xymatrix{{B^3} \ar[r]^h & {A}}} & {\xymatrix{{B^2} \ar[r]^{\varphi} & {A}}} &
{\xymatrix{{B}
\ar[r]^{u}& {A}}}
\end{array}
\ee
The pair ($g,h$) satisfies equations (\ref{bicocyc})-(\ref{bicompat}), $\varphi$ and
$u$ satisfy the conditions (\ref{phitriv}) and (\ref{phisym}), together with the additional
conditions
\be
\label{ctriv} 
c(b,b') = \frac{u (b + b')}{u (b) \: u (b')}
\ee
and 
\be
\label{u2}
u(b)^2 = \varphi (b, b)
\ee 
where $c$ is defined, in terms of the triple $(g,h,\varphi)$, by the
equations (\ref{lambdabiex}) and (\ref{defc}). 
A trivialization of the alternating biextension defined by  $(g, h, \varphi, u)$ is described
by a map  $k:\nolinebreak B^2 \la A$ which satisfies the equations (\ref{bitriv}) and (\ref{ktriv}),
and the additional condition
\be
\label{ktriv2}
\begin{array}{rcl}
k(x,y) \: k(y,x) & = & \varphi (x,y)  \\
\end{array}
\ee

\end{proposition}

\begin{remark} {\rm
Here is the connection between this second cocyclic description of an alternating biextension and
the original one in terms of a triple ($g,h,u$), where $(g,h)$ again  satisfies the biextension
cocycle conditions (\ref{bicocyc})-(\ref{bicompat}) and $u$ satisfies equations (\ref{symcoc1}) and
(\ref{symcoc2}). Starting from such a triple  ($g, h, u$), one defines an
anti-symmetry isomorphism map
$\varphi: B^2 \la A$ by 
\begin{eqnarray}
\label{newphi}
\varphi (b,b') & = & \frac {u(b+b')}{u(b) \: u(b') \: \lambda (b,b')}
 \end{eqnarray}
It is immediate that
$\varphi(b,b')$ also satisfies  the symmetry condition (\ref{phisym}), and equations (\ref{phitriv})
are consequences of (\ref{symcoc2}). In fact, the triple $(g,h,\varphi)$  describes in
cocyclic  terms the anti-symmetric biextension determined by the alternating biextension
$(g,h, u)$. Furthermore, by definition of
$\varphi (b,b')$,  the cocycle
$c(b,b')$ defined by (\ref{defc})  satisfies the trivialization condition (\ref{ctriv}).  
In order to verify that the quadruple
$(g, h,
\varphi, u)$  associated to the triple ($g, h, u$) and to the map $\varphi$
(\ref{newphi}) satisfies the conditions of proposition  \ref{quadrup}, it  suffices to
check  condition (\ref{u2}), in other words that the equation
\begin{eqnarray}
\label{u22}
 u(b)^2 & = & \frac{u(2b)}{u(b)^2 \: \lambda (b,b)} 
\end{eqnarray}
 is satisfied. By specialization to the case $b_2 = -b_1, \: b_3=b_1$, (\ref{symcoc2}) yields the
equation
\begin{eqnarray*} \frac{u(b)^3 u(-b)}{u(2b)} &  =  & \frac{g(b, -b ; b) \: h(b; b, -b)}{\lambda
(b,b) \:
\lambda (-b,b)}
\end{eqnarray*}
 Substituting in this equation  the values for $u(-b)$ and for $\lambda$ given by (\ref{symcoc1})
and (\ref{lambdabiex}) yields the requisite formula (\ref{u22}).
}
\end{remark}

\section{The commutator of $\cc$ as a weak biextension}
\label{sec:comweak}

\setcounter{equation}{0}%

\hspace{.7cm}We are now ready to describe 
 the universal coefficient exact sequence (\ref{uct3}) in geometric terms. 
Let $\cc$ be a {\it gr}-category, as defined in \S \ref{sub:uct},   with invariants $B$ and $A$
satisfying the hypothesis \ref{hyp1}. Suppose  that we have chosen for each $x \in B$, as we did
above,  a representative object
$X_x$ of
 $\cc$ in the
isomorphism class of $x$. To $\cc$ we  associate
the $A$-torsor $E$ above  $B \times
B$, which Deligne \cite{pd:symbole} calls the commutator of
$\cc$, whose fibre  above $(x,y) \in B^2$ is the set
\begin{equation}
\label{biexdef}
E_{x,y} = 
Isom_{\,\cc}(X_y X_x, X_x X_y) 
\label{pr:Ex,y}
\end{equation}
of arrows from $X_y X_x$ to $X_x X_y$. Composing the elements of $E_{x,y}$ on the right
with automorphisms of $X_y X_x$, 
viewed as elements of $A$, makes
$E$ into a right $A$-torsor on $B \times
B$. Alternate choices for the
 representative objects $X'_x$ and
$X'_y$ of $x$ and $y$ yield an $A$-torsor $E'$ on $B \times
B$
isomorphic to $E$.

\bigskip

The main result of this section is the following proposition.

\begin{proposition}
\label{pr:Eweak}
The $A$-torsor $E$ associated 
 to the monoidal category $\cc$ is
endowed with a natural structure of a weak
biextension of $B \times B$ by $A$.
\end{proposition}
{\bf Proof:}  In order to simplify the notation, we will in the following discussion denote by $X,
Y, Z$,
 {\it etc ...} the chosen
representatives $X_x, X_y, X_z$ in $\cc$ of the elements $x,y, z \in B$. Let $u: YX
\longrightarrow XY$ and $v: Y X'
\longrightarrow X' Y$ be given
elements in $E_{x,y}$ and $E_{x',y}$. Their
partial sum $u
\stackrel{1}{+} v$ 
(\ref{eq:+1}) is defined to be the section of $E_{xx',y}$ determined as follows. Consider 
the following composite arrow, in
which the unlabelled arrows are the associativity
isomorphisms.

\begin{equation}
\label{u+1v} 
Y (X  X')
\rightarrow (Y
 X) X'
\stackrel{u}{\rightarrow} (X
Y) X' \rightarrow X
(Y X')
\stackrel{v}{\rightarrow} X
 (X' Y) \rightarrow (X
 X')  Y
\end{equation}
Reverting temporarily to our standard notation for objects, this is the middle arrow in the
composite map 
\be
\label{u+1v2}
X_y X_{xx'} \la X_y (X_x X_{x'}) \la (X_x X_{x'}) X_y \la X_{xx'} X_y
\ee
whose other arrows are defined by (\ref{cxy}).
If $w: Y'X \longrightarrow XY'$ is
another arrow in $\cc$,  the partial
sum $u \stackrel{2}{+} w$ (\ref{eq:+2}) is the
composite map constructed in a similar way from the arrow
\begin{equation}
\label{u+2w}
(YY')X \rightarrow Y(Y'X)
\stackrel{w}{\rightarrow} Y(XY') \rightarrow
(YX)Y' \stackrel{u}{\rightarrow} (XY)Y'
\rightarrow X(YY').
\end{equation}
Observe that the composite arrows (\ref{u+1v}) and
(\ref{u+2w}), which are built out of intertwining
associativity and commutativity isomorphisms,
are  the familiar boundary arrows in the two hexagons occurring as axioms for braided monoidal
categories\footnote{ We do not, however, assume that our monoidal category
$\cc$ is braided.} \cite{js:braided}. It is therefore not surprising that the diagrams describing the
required associativity and compatibility conditions
 for our
composition laws
 $\stackrel{1}{+}$ and  $\stackrel{2}{+}$ are closely related to
some of the higher braiding axioms embodied in
the definition of a braided 2-category
\cite{kv:2-cat}.  Specifically,  for
each set of elements $ a:WX 
\longrightarrow XW , \:  
b:WY \longrightarrow YW$ and $ 
c:WZ \longrightarrow ZW$  of
$E_{x,w}$, $E_{y,w}$ and $E_{z,w}$, we must
consider a non-strict version of the tetrahedral diagram analogous
to the diagram of type 
$(\bullet \otimes (\bullet
\otimes \bullet \otimes \bullet))$
associated in
\cite{kv:2-cat} \S 6 to the objects $W, X, Y, Z$ of
$\cc$ (see also
\cite{ba-ne:higher}).
 In order to take into account 
the associativities, this requires that we double certain edges of this
diagram (in fact precisely those edges which are thickened in  the
diagram appearing in \cite{kv:2-cat} \S 6 ).
\begin{equation}
\xymatrix{&& *+<2ex>{XYZW} &&\\
&&&&\\
{WXYZ} \ar[uurr] \ar@<1ex>[uurr] \ar@{-->}[rrrr]
\ar@<-0.5ex>[drr]^{a(YZ)}\ar@<-1.5ex>[drr]_{(aY)Z} &&&&{XYWZ} 
\ar@<-0.5ex>[uull]^{X(Yc)} \ar@<-1.5ex>[uull]_{(XY)c}\\
&&{XWYZ}  
\ar[uuu]
\ar@<-0.5ex>[urr]^{X(bZ)}
\ar[uuu]\ar@<-1.5ex>[urr]_{(Xb)Z}&&}
\label {pic:assoc1}
\end{equation}
\vskip0.5cm
\noindent There are now five edges incident with each
vertex. Replacing each of these vertices by
the corresponding commutative pentagon in $\cc$,  we may now attach one of
the incident edges to
each of the five vertices of each pentagon.
Taking into account the labels given to
certain arrows, this can be done in a unique
manner, if we require that exactly three edges be
incident to each vertex of each pentagons.In our context, this diagram has the following
interpretation. The axioms on $\cc$  ensure that every face of the polyhedron is now 
 commutative, except possibly for the face 
comprising the two arrows between  $WXYZ$ and
$XYZW$. The latter face is the following
square, whose vertical arrows are the
associativity isomorphisms:

\[
\xymatrix{
{W(X(YZ))} \ar[dd] \ar[rr]^{a \stackrel{1}{+} (b \stackrel{1}{+} c)}
&&  {(X((YZ))W)} \ar[dd]\\
\\
{W((XY)Z)} \ar[rr]^{(a \stackrel{1}{+} b) \stackrel{1}{+} c}
&& {((XY)Z)W)}} 
\]
 Since all other faces of diagram (\ref
{pic:assoc1}) commute, so does this face. This finishes the proof
that the composition law
$\stackrel{1}{+}$ is indeed associative. The
associativity of  the law $\stackrel{2}{+}$ is
obtained in a similar manner, by starting
instead from the non-strict version of diagram 
$ ((\bullet
\otimes \bullet
\otimes
\bullet) \otimes \bullet)$ of
\cite{kv:2-cat} \S 6.

\medskip

The compatibility between the composition laws
$\stackrel{1}{+}$ and $\stackrel{2}{+}$ is proved
in a similar manner, starting instead from a non-strict version\footnote{I owe to E. Getzler the observation that  the non-strict versions
of  diagrams 
$(\bullet \otimes (\bullet
\otimes \bullet \otimes \bullet))$ and  $ ((\bullet
\otimes
\bullet) \otimes (\bullet \otimes \bullet))$ appear respectively as figures 4 and 5 of \cite{dbn}} of 
 diagram  $ ((\bullet
\otimes
\bullet) \otimes (\bullet \otimes \bullet))$
of \cite{kv:2-cat} \S 6.
This is the diagram

\begin{equation}
\label{pic:compat1}
\xymatrix{
{XYZW} \ar@<-1ex>[dd]^{(Xd)W}
\ar@<-2ex>[dd]_{X(dW)}
\ar[dddrrr]
\ar@{-->}[drrr]
\ar@<-0.5ex>[rrrrrr] \ar@{-->}@<0.5ex>[rrrrrr]&&&&&& ZWXY \\
&&& {XZWY} \ar@{-->}[urrr] \ar@{-->}[drrr]^(.3){(cW)Y} 
\ar@{-->}@<-1ex>[drrr]_(.3){c(WY)}&&&\\{XZYW}\ar@{-->}[urrr]^(.7){X(Zb)} 
\ar@{-->}@<-1ex>[urrr]_(.7){(XZ)b}\ar@<-0.5ex>[drrr]^{(cY)W}
\ar@<-1.5ex>[drrr]_{c(YW)}&
&
&
&
&
& {ZXWY}\ar[uu]^{(Za)Y}  \ar@<-1ex>[uu]_{Z(aY)}\\
&
&
&
{ZXYW}
\ar[urrr]^{(ZX)b}
\ar@<-1ex>[urrr]_{Z(Xb)}
\ar[uuurrr]
&
&&}
\end{equation}
associated to four given arrows $a:XW \longrightarrow WX , \: 
b:YW \longrightarrow WY , \:
c:XZ \longrightarrow ZX , \: 
$ and $d:YZ \longrightarrow ZY$.
Every vertex of diagram  $ ((\bullet
\otimes
\bullet) \otimes (\bullet \otimes \bullet))$
of \cite{kv:2-cat} \S 6 has now been replaced
by the corresponding associativity pentagon. The only {\it a priori} non
commutative part of our diagram is the square involving the top two
horizontal arrows between $XYZW$ and $ZWXY$ (and appropriate
associativity arrows).  In particular, the following
commutative triangles which we extract from diagram (\ref {pic:compat1})
yield for us the appropriate labels for these horizontal arrows

\begin{equation}
\begin{array}{cc}
\mbox{\xymatrix{
{XYZW} 
\ar[rr]^{(a \stackrel{1}{+} b) \stackrel{2}{+} (c\stackrel{1}{+} d) }
\ar[ddr]_{(c
\stackrel{1}{+} d)W} && {ZWXY} \\ &&&
\\&{ZXYW}
\ar[uur]_{(a \stackrel{1}{+} b)W}&}}
 & 
\mbox{\xymatrix{
{XYZW} \ar@{-->}[rr]^{ (a \stackrel{2}{+} c) \stackrel{1}{+}(b
\stackrel{2}{+} d)}
\ar@{-->}[ddr]_{X(b
\stackrel{2}{+} d)}&& {ZWXY} \\
&&&\\
& {XZWY} \ar@{-->}[uur]_{(a \stackrel{2}{+} c)Y}&}}
\end{array}
\end{equation}

\noindent The square involving the top two horizontal
arrows of (\ref {pic:compat1}) may therefore be portrayed, with
certain associativity morphisms neglected, as 

\begin{equation}
\xymatrix{
{XYZW} \ar[rrr]^{(a \stackrel{1}{+} b) \stackrel{2}{+}
(c\stackrel{1}{+} d)} \ar[d] &&& {ZWXY} \ar[d]\\
{XYZW} \ar@{-->}[rrr]^{(a
\stackrel{2}{+} c) \stackrel{1}{+}(b
\stackrel{2}{+} d)}&&& 
{ZWXY}}
\end{equation}
This commutes, since all the other faces of diagram (\ref
{pic:compat1}) do. This proves that the partial group laws
$\stackrel{1}{+}$ and $\stackrel{2}{+}$ in $E$  respectively defined by
(\ref {u+1v}) and (\ref {u+2w}) are   compatible and therefore finishes
the proof of proposition \ref {pr:Eweak}.

\begin{remark} 
\label{cocyclbiext}
{\rm The chosen  arrows $c_{x,y}$ (\ref{cxy}) determines a section 
\be
\label{cxy1}
d_{x,y}: X_y X_x \la X_{yx} = X_{xy} \la X_x X_y
\ee
of the torsor  underlying the commutator $E$ (\ref{pr:Ex,y}) of $\cc$. The cocycles which express
the partial group laws $\stackrel{1}{+}$ and $\stackrel{2}{+}$ of $E$ in terms of this section may
now be made explicit. Let the sections $u$ and $v$ of $E_{x,y}$ and
$E_{x',y}$ be the chosen morphisms
$d_{x,y}$ and
$d_{x',y}$. 
 The following  commutative diagram, in which the  unlabelled  horizontal arrows are all  identity
maps,  expresses the  composite map (\ref{u+1v2}) in terms of automorphism of $X = X_{xx'y}$.
\begin{equation}
\label{defcoc}
\hspace{-1cm}{
\xymatrix{
{X_y(X_xX_{x'}) } \ar@{-}[d] \ar[r]^{a_{y,x,x'}} & {(X_yX_x)X_{x'}} \ar@{-}[d]
\ar[r]^{d_{x,y}X_{x'}}  & {(X_xX_y)X_{x'}} \ar@{-}[d] \ar[r]^{a^{-1}_{x,y,x'}}& {X_x(X_yX_{x'}) }
\ar@{-}[d]
\ar[r]^{X_xd_{x',y}} & {X_x(X_{x'}X_y)} \ar@{-}[d] \ar[r]^{a_{x,x',y}} & {(X_xX_{x'})X_y} \ar@{-}[d]
\\ {X_y X_{xx'}} \ar@{-}[d] & {X_{yx}X_{x'}} \ar@{-}[d] \ar@{-}[r] & {X_{xy}X_{x'}} \ar@{-}[d]
& {X_xX_{yx'}} \ar@{-}[d]\ar@{-}[r] & {X_xX_{x'y}} \ar@{-}[d] & {X_{xx'}X_y} \ar@{-}[d]\\
{X}\ar[r]_{f_{y,x,x'}}
&{X}\ar[r] 
&{X} 
\ar[r]_{f^{-1}_{x,y,x'}}&{X}\ar[r]
&
{X}\ar[r]_{f_{x,x',y}}
&
{X}}}
\end{equation}
The cocycle $g(x,x' \, ; y)$ which describes the
partial sum
$\stackrel{1}{+}$ may now be read off from the lower horizontal map of this diagram as the map
$g:B^3 \la A$  obtained 
 from  the 3-cocycle $f(x,y,z)$ by shuffling $y$ through $x,x'$, in other words by the formula
\be
\label{biext:1}
g(x,x' \,;  y)  = \frac{f(x,x',y) \: f(y,x,x')}{f(x,y,x')}
\ee
Starting instead  from the definition (\ref{u+2w}) for $\stackrel{2}{+}$, one sees that the second
cocycle
$h$ occurring in the definition of the commutator biextension is obtained by shuffling $x$ through
$y,y'$ in the opposite direction, in other words by 
the rule 
\be
\label{biext:2}
h(x; \: y,y') = \frac{f(y,x,y')}{f(x,y,y') \: f(y,y', x)}
\ee
 That the pair ($g(x, y \, ;z)\, , \, h(x; \, y, z)$) satisfies the cocycle conditions
for a weak  biextension
 follows from the previous discussion. 
It could also have been be proved directly by repeated use of the 3-cocycle condition for $f$.}
\end{remark}

\section{The trilinear map associated to a monoidal category}
\label{par: trilinmonoid}

\setcounter{equation}{0}%

\hspace{.7cm}We now examine  the conditions under which the commutator weak biextension
(\ref {pr:Ex,y}) is a genuine biextension. In view of
proposition \ref {pr:Eweak}, it suffices to check that  both
partial multiplication laws 
on $E$ are 
commutative.
 In contrast to
the associativity and compatibility conditions, the commutativity conditions are
not automatically satisfied. At the cocycle level, it is immediate that each of the two
commutativity axioms  leads to  the  following condition on $f$
 \[ \frac{f(x,y,z) \: f(z,x,y) \: f(y,z,x)}{f(x,z,y) \: f(z,y,x) \: f(y,x,z)} = 1
\]
 The expression $\varphi(x,y,z)$ defined by the left-hand side of this 
equation is simply the evaluation of the 3-cocycle $f$ on the triple Pontrjagin product cycle
$x.y.z \in H_3(B)$ of classes $x,y,z \in H_1(B)$. It is well-known that $\varphi$ is a
trilinear alternating map. This proves the following proposition

\begin{proposition}
\label{pr:trilin}
For any pair of   abelian groups  $A$ and $B$, the weak commutator
$A$-biextension
$E$ on
$B
\times B$ (\ref {pr:Ex,y}) associated  to a monoidal
category
$\cc$ satisfying hypothesis \ref {hyp1} described by a $3$-cocycle $f(x,y,z)$ is  a genuine 
biextension of
$B
\times B$ by $A$ if and only if the alternating map

\begin{equation}
\label{pic:wedgetrilin}
\xymatrix{
{B \wedge B \wedge B } \ar[r]^<<<<<{\varphi_{\cc}} & {A}
}
\end{equation}
defined by
\be
\label{eq:trilin}
\varphi_{\cc} (x, y,z) = \frac{f(x,y,z) \: f(z,x,y) \: f(y,z,x)}{f(x,z,y) \: f(z,y,x) \: f(y,x,z)}
\ee
is trivial.
\end{proposition}

This statement may also be understood in geometric terms, without appealing to the preferred
section
$d$ (\ref {cxy1}) of $E$, by considering the following diagram, which is built by pasting together
two diagrams of type (\ref{defcoc}) associated to  the partial sums $u \stackrel{1}{+} v$ and  
$v \stackrel{1}{+} u$ of a pair of arbitrary  summable sections $u,v$ of $E$. We no longer
assume here that the vertical arrows are the specific morphisms $d_{x,y}$ and $d_{x',y}$, so that
the automorphisms of $X$ which $u$ and $v$  determine are not necessarily the identity arrows.
The latter are denoted $g_u$ and $g_v$.

\begin{equation}
\label{pic:com1}
\hspace{-1cm}
\xymatrix@=4ex{
{X_{y}(X_{x'}X_{x})}
\ar@{-}[ddd]
\ar[r]^{a_{y,x',x}} & {(X_{y}X_{x'})X_{x}}
\ar[r]^{vX_x} & {(X_{x'}X_{y})X_{x}}
\ar[r]^{a^{-1}_{x',y,x}} & {X_{x'}(X_{y}X_{x})}
\ar[r]^{X_{x'}u}& {X_{x'}(X_{x}X_{y})}
\ar[r]^{a_{x',x,y}} & {(X_{x'}X_{x})X_{y}}
\ar@{-}[ddd]
\\  &  \ar@{-}[u] {X_{yx'}X_x}
\ar[r]^{g_v X_x}& {X_{x'y}X_x}
\ar@{-}[u]& {X_{x'} X_{yx}} \ar@{-}[u]
\ar[r]^{X_{x'}g_u} &  {X_{x'} X_{xy}} \ar@{-}[u]& \\
& {X} \ar@{-}[u]
\ar[r]^{g_v} & {X} \ar@{-}[u] 
\ar[r]^{f_{x',y,x}^{-1}}
& {X}  \ar@{-}[u] \ar[r]^{{}^{x'}\hspace{-4pt}g_u}& {X} \ar@{-}[u]
\ar[d]^{f_{x',x,y}}& \\
{X_yX_{x'x}} \ar@{-}[ddd] 
\ar[r] & {X} \ar[u]^{f_{y,x',x}}
\ar[d]_{f_{y,x,x'}}&&& {X} \ar@{-}[r] 
&
{X_{x'x}X_{y}} \ar@{-}[ddd]
\\ & {X} \ar@{-}[d]
\ar[r]_{g_u}& {X} \ar@{-}[d]
\ar[r]_{f_{x,y,x'}^{-1}}& {X} \ar@{-}[d] \ar[r]_{{}^{x}g_v}& {X}
\ar@{-}[d]
\ar[u]_{f_{x,x',y}}&\\
 & {X_{yx}X_{x'}}
\ar@{-}[d]
\ar[r]_{g_u X_{x'}} & {X_{xy}X_{x'}} \ar@{-}[d]&
{X_x X_{yx'}} \ar@{-}[d] \ar[r]_{X_{x}g_v}
& {X_x X_{x'y}} \ar@{-}[d]& \\
{X_{y}(X_{x}X_{x'})}
\ar[r]_{a_{y,x,x'}} & {(X_{y}X_{x})X_{x'}}
\ar[r]_{uX_{x'}} & {(X_{x}X_{y})X_{x'}}
\ar[r]_{a_{x,y,x'}^{-1}}& {X_{x}(X_{y}X_{x'})}
\ar[r]_{X_{x}v}
&{X_{x}(X_{x'}X_{y})}\ar[r]_{a_{x,x',y}}
&{(X_{x}X_{x'})X_{y}}
\\
&&&&&\\}
\end{equation}
 
\vskip 0.5cm
\noindent All cells in this diagram are commutative,
 except possibly
the large inner one, composed of automorphisms of the object $X = X_{xx'y}$. Since the group
$A = Aut(X)$ is abelian, all the arrows lying on the boundary of this inner cell may be freely moved
past one another. Furthermore, under hypothesis
\ref{hyp1}, the elements ${}^{x'}\!g_u$ and $g_u$ in the inner section
coincide, so that they cancel each other out, and
similar cancellation occurs between
${}^{x'}\!g_v$ and 
$g_v$. Keeping track of the orientations of the arrows, the
obstruction to the commutativity of the inner region boils down to the
triviality of the expected element 
$\varphi_{\cc} (x, x',y)$ (\ref{eq:trilin})
of $A$.

\bigskip

\begin{remark}
{\rm 
 The specific arrows $u$ and $v$ chosen in the
construction of diagram (\ref{pic:com1}) played no role in the definition
of the map  $\varphi_{\cc}$
(\ref {pic:wedgetrilin}), which only depended on the $3-$cocycles
$f_{x,x',y}$ determined by the associativity data $a_{x,x',y}$ in $\cc$. 
Another set of choices for the vertical maps $X_xX_y
\la X_{xy}$ will yield a $3-$cocycle $f'$ cohomologous to $f$, and which
therefore leaves unchanged  the induced map (\ref {eq:trilin}). The
axioms for $\cc$ also ensure that alternate choices for the
representative objects $X_x$ of $x \in B$ have no effect in this construction.}
\end{remark}

\section{The  alternating structure on the commutator biextension}
\label{par: altcom}

\setcounter{equation}{0}%

\hspace{.7cm}There remains yet one additional element of structure of $E$ to be made explicit. We
now assume that the map (\ref {pic:wedgetrilin}) associated to the given
monoidal category $\cc$ is trivial, so that by proposition 
\ref {pr:trilin} the commutator torsor $E$ (\ref {pr:Ex,y}) is a
genuine biextension. 

\begin{proposition}
\label{pr:altern}
Let  $\cc$ be a monoidal category with invariants the abelian groups
$\pi_{0}(\cc) = B$ and  $\pi_{1}(\cc , I) = A$,  satifying hypothesis
\ref{hyp1}, and whose associated trilinear map (\ref{eq:trilin}) is
trivial. The associated commutator biextensions $E$ of $B \times B$
by $A$ is alternating.
\end{proposition}

\noindent{\bf Proof:}  We give here  a geometric proof of this assertion, which one was in any case
 led to expect by the discussion in
\S
\ref{sub:uct}. We begin with a geometric proof of the following weaker assertion, which can also
be deduced  from the cocyclic description (\ref{biext:1})-(\ref{biext:2}) of $E$.

\begin{lemma}
\label{lem:antisym}
The commutator biextension $E$ of $\cc$ is anti-symmetric.
\end{lemma}
\noindent {\bf Proof:} Consider the section $\sigma$ of  $E \wedge s^{*}E$ above $B
\times
B$ (where $s: B\times B \longrightarrow B \times B$ is the map which
permutes the factors) defined by the rule  which assigns to any pair of elements
$x,y
\in B\times B$ the  element 
\[\sigma(x,y) = t \wedge v
\] in
$ \mathrm{Isom}(X_y X_x, X_x, X_y) \wedge \mathrm{Isom}(X_x X_y, X_y X_x) 
$
where $v:X_x X_y \longrightarrow  X_y X_x$ is the inverse of
the arrow $t: X_y X_x \longrightarrow X_x, X_y$. This section of $E \wedge
s^{*}E$ does not in fact depend on the choice of a specific map $t: X_y X_x \longrightarrow X_x,
X_y$.
In order to check that  $\sigma$ trivializes $E \wedge
s^{*}E$  as a biextension, it must be verified that it
 is multiplicative in
each of its two variables.  In the first variable, this
boils down to the obvious assertion that for a given pair of sections 
$t:X_y X_x
\longrightarrow X_x X_y$ and   $ t':X_y X_{x'} 
\longrightarrow X_{x'} X_y$ of $E_{x,y}$ and $E_{x',y}$ with
inverses $v$ and $v'$, the inverse in $\cc$ of the composite map (\ref{u+1v})
$t \stackrel{1}{+} t'$
\begin{equation}
\xymatrix{
{X_y X_{x'}X_x}\ar[rr]^{t}
&
&{X_xX_yX_{x'}}
\ar[rr]^{t'}
&
&{X_xX_{x'}X_y}}
\end{equation}
is the map  $v \stackrel{2}{+} v'$ (\ref{u+2w})
\begin{equation}
\xymatrix{
{X_xX_{x'}X_y}\ar[rr]^{v'}
&&{X_x X_y X_{x'}}
\ar[rr]^{v}
&
&
{X_yX_xX_{x'}}}
\end{equation}
The multiplicativity of $\sigma$ in the second variable is verified in a similar manner, so that the
lemma is proved.

\bigskip 

Since $E$ is anti-symmetric, its restriction $\Delta E$ to the diagonal is a (commutative)
extension of $B$ by $A$, for which the group law $\star$ may be described explicitly by the
following rule. Let $u$ and $v$ respectively be sections of $E_{x,x}$ and $E_{y,y}$, in other
words  arrows
\be
\begin{array}{ccc}
\mbox{$u: X_x X_x \la X_x X_x$} &,& \mbox{$v: X_y X_y  \la X_y X_y$}
\end{array}
\ee
in $\cc$. The arrow $u \star v : X_{xy} X_{xy} \la  X_{xy}   X_{xy}$
is determined by the composite map
\be
\label{ustarv}
\xymatrix{
{X_x X_y X_x X_y } \ar[rr]^{X_x \alpha X_y} & &{(X_xX_x)(X_yX_y)} \ar[r]^{uv}
 &{(X_xX_x)(X_yX_y)}\ar[rr]^<<<<<<<<<{X_x \alpha^{-1} X_y} & & {X_x X_y X_x X_y} }
\ee
for some arbitrarily chosen  arrow $\alpha :X_y X_x \la X_x X_y$. We know that $(\Delta E)^2$,
being the restriction to the diagonal of
$E
\wedge s^{*}E$, is trivialized by the restriction $\Delta \sigma$ of  $\sigma$ to the diagonal. For
any
$x
\in B$, it therefore follows that  $\sigma(x)$ is  the element $t \wedge t^{-1}$ in the set 
$\mathrm{Isom}(X_xX_x,X_xX_x)
\wedge
\mathrm{Isom}(X_xX_x,X_xX_x) $,  for some arbitrarily chosen arrow
$t: X_x X_x
\la X_x X_x$. Choosing for $t$ the identity self-arrow
$1_{X_xX_x}$,  we may therefore set $\sigma (x) =1_{X_xX_x}
\wedge 1_{X_xX_x}$. Consider now the section $\tau$ of $\Delta E$ defined by setting 
\[ \tau (x) =  1_{X_xX_x} \]
The equation $\tau (x) \star \tau (y) = \tau (x+y)$ follows from the definitions, so
that the section
$\tau$ splits
$\Delta E$ as an extension. The formula 
\[ \sigma (x) = \tau (x) \wedge \tau (x) \]  
is also immediate.  By proposition
\ref{altantisym}, the section
 $\tau$ of $\Delta E$  therefore induces an alternating biextension structure on the anti-symmetric
biextension
$E$.
\bigskip 

It is easily verified that a trivialization of $E$ compatible with its alternating
biextension structure determines a strict Picard  structure on the monoidal
category  $\cc$. As observed in the introduction, it follows from \cite{pd:sga4} \S 1.4 that such
strict Picard categories are always trivial, since they are classified up to equivalence by the
(trivial) group $\mathrm{Ext}^2(B, A)$. This may be spellt out as follows.
\begin{corollary}
\label{cor:alt}
 Let $\cc$ be a monoidal category satisfying the conditions of
proposition (\ref{pr:altern}) whose associated biextension $E$ is  
trivial (as an alternating biextension). Then the monoidal structure on 
$\cc$ is trivial.
\end{corollary}

In more concrete terms,  consider a monoidal category $\cc$ determined by a crossed module
$\delta:M \la N$, for which
$\ker(\delta) = A$ and $\mathrm{coker}(\delta) = B$, and for which the $B$-module structure on $A$
is trivial. The fibre
$E_{x,y}$ of the commutator biextension $E_{\cc}$ above $(x,y) \in B^2$ is the set
$\delta^{-1}(k_{x,y})$, where $k_{x,y}
\in K =
\mathrm{im}(\delta)$ is the image of $(x,y)$ under the commutator map associated to the central
extension
\[
\mathrm{0} \la K \la N \la B \la \mathrm{0}
\]
 A trivialization of $E$ determines in  a map 
\be
\label{crmodule}
 \{ \: , \: \}: N \times N \la M
\ee
by associating to a pair of elements $X$ and $Y$ in $N$ with projection $x$
and
$y$ in
$B$  the
 element $m \in  \delta^{-1}(k_{x,y}) \subset M$ determined by the trivialization. When the
trivialization of $E$ is compatible with the  anti-symmetric biextension structure on $E$, the map
\ref{crmodule} determines a stable crossed module structure \cite{mc.gen} on $\delta:M \la N$.
Compatibility of the trivialization of $E$ with the alternating structure on $E$ yields the additional
relation $\{n,n\} = 0$ for all $n \in N$. In the present context, Deligne's result \cite{pd:sga4}
asserts that the given crossed module is then equivalent to the crossed module $P \la Q$ determined
by a (splittable) exact sequence of abelian groups
\[ 
\mathrm{0} \la A \la P \la Q \la B \la 0
\]

\bigskip

We end this section with a brief discussion in cocyclic terms of proposition \ref{pr:altern}
and of its  corollary. When the cocycle pair $(g,h)$  has been defined in terms of a
3-cocycle
$f$ by equations (\ref{biext:1}) and (\ref{biext:2}), both the map $\lambda (b_1, b_2)$
(\ref{lambdabiex}) and the right-hand terms of equation (\ref{phitriv}) are trivial.
Taking into account (\ref{defc}),  it follows  that the quadruple
$(g,h, 1, 1)$ satisfies the conditions of proposition \ref{quadrup} so that it defines an
alternating biextension structure on  the biextension $(g,h)$. The underlying triple $(g,h,1)$
then satisfies the equivalent conditions (\ref{symcoc1}) and (\ref{symcoc2}), as asserted in
proposition
\ref{pr:altern}. These two conditions on the triple $(g,h,1)$ may also be verified directly 
without introducing explicitly the full quadruple $(g,h,1,1)$.

\bigskip

A discussion in similar
terms of corollary
\ref{cor:alt} goes as follows. Suppose that the quadruple $(g,h,1,1)$ associated to a cocycle $f$
is trivial, so that there exists a  map $k(b_1, b_2)$ satisfying conditions
(\ref{bitriv}), (\ref{ktriv}) and (\ref{ktriv2}). 
In that case the pair ($f,k$) defines an element in the trivial group
$\mathrm{Ext}^2(B, A)$ so that, as observed at the end of remark \ref{rem13}, there exists
an
$A$-valued   2-cochain
$l(b_1,b_2)$ on
$B$ for which 
\begin{eqnarray}
k(b_1, b_2) = \frac{l(b_1,\, b_2)}{l(b_2,\, b_1)} 
\end{eqnarray}
and such that the 3-cocycle $f$ is the coboundary of $l$. The latter assertion is the content of
corollary
\ref{cor:alt}. We note in passing that the category $\cc$ described by the 3-cocycle $f$ is
braided if and only if   there exists a map $k$ which trivializes the pair  $(g,h)$ as a
biextension (without taking into account the alternating structure), in other words which 
satisfies equations (\ref{bitriv}) but not (\ref{ktriv}).

\section{From monoidal categories to monoidal stacks}

\setcounter{equation}{0}%

\hspace{.7cm}The previous analysis of monoidal categories {\it via} the universal coefficient theorem
 extends to a classification of group-like monoidal stacks in groupoids
 (also called {\it gr}-stacks)
in a general topos $T$,
 as discussed
in \cite{lb:2-gerbes} \S 7, to which we refer for the requisite definitions. One is given a pair of
 abelian groups $B$ and $A$ of $T$, with $A$ viewed as
a trivial $B$-module. The discussion in \S \ref{sub:uct} generalizes to the assertion that monoidal
stacks $\cc$ of $T$ with invariants $\pi_0(\cc)$ and $\pi_1(\cc)$  respectively isomorphic to
$B$ and $A$ and satisfying hypothesis \ref{hyp1} are classified by the hypercohomology group
$H^3(B, A)$.  The difference between this hypercohomology group and the ordinary
cohomology group
$H^3(B,A)$ is analyzed by the first quadrant spectral sequence
\be
\label{hypcoh1}
E^{p,q}_1 = H^q(X_p, A) \Longrightarrow H^{p+q}(B, A)
\ee
whose initial term is the $A$-valued cohomology of the degree $p$ component $B^p$ of the classifying
space $X_\ast$ of $B$. From the geometric point of view which concerns us here, the distinction between
the hypercohomology group  and the naive cohomology group $H^3(B,A)$ built from
cochains $B^3 \la A$, and  which classify {\it gr}-categories with invariants 
$B$ and $A$, is reflected in the two sets of obstructions whose vanishing is necessary
in order to carry out in the stack case the construction
of the 3-cocycle associated to the monoidal category $\cc$. 
 The first  of these arises when one attempts to choose,
 for each section $x$ of $B$  above some object $S$ of $T$,  an object $X_x$ in the 
fiber category $C_S$ of $\cc$ above $S$. Suppose that this obstruction vanishes, 
so that  the requisite objects $X_x$ exist for all
sections $x$ of $B$. Just as in (\ref{cxy}),  one  may  then attempt to choose, for every pair of
sections $x,y$
 of $B$ above a given object 
$S$ of $T$, a morphism $c_{x,y}: X_x X_y \la X_{xy}$ in $\cc_S$. If the obstruction to achieving
 this
also vanishes, then one can construct as in (\ref{diag:f}), an $A$-valued 3-cocycle and therefore
classify  by the 
naive group $H^3(B,A)$ the stacks  $\cc$ for which both sets of obstructions vanish. These
obstructions  do not however vanish in general,  but the
 stack axioms, and the definitions of the objects
$\pi_i(\cc)$ in $T$ ensure nevertheless that they both vanish locally\footnote{in
other words after base change.}. The invariants
which describe them are therefore of a cohomological nature. We refer to  \cite{lb:2-gerbes} \S 7 for 
a further discussion of these invariants, and simply observe here that they live  respectively in the
terms  $E^{1,2}_1$ and $E^{2,1}_1$  of the spectral sequence (\ref{hypcoh1}), while the naive
 cohomology group
$H^3(B,A)$ is its  $E^{3,0}_2$ term. The comparison
 between the naive
cohomology group and the associated hypercohomology group (in other words  between
the classification of {\it gr}-stacks and that of {\it gr}-categories) 
  thus boils down to the analysis of the edge-homomorphism map  $E^{3,0}_2 
\la H^3$ in the spectral sequence (\ref{hypcoh1}).

\medskip
Another change occurs when one passes from categories to stacks.
While the homology of the abelian group $B$ of a topos $T$ is 
given by the same formulas as for an abstract group,
 the relation between the homology
and the hypercohomology of  $B$ is now more complicated, since the 
universal coefficient theorem must now be replaced by  the universal coefficient spectral sequence
\be
\label{hypcoh2}
E^{p,q}_2 = \mathrm{Ext}^p(H_q(B), A) \Longrightarrow  H^{p+q}(B, A)
\ee
In the abstract group situation, this spectral sequence reduces to the ordinary
 universal coefficient theorem, since in the category of 
abelian groups the groups $\mathrm{Ext}^p$ vanish whenever $p > 1$ so that the spectral 
sequence degenerates.  In a topos,  no essential change occurs at the level of degree 2 cohomology,
so that the analysis of central extensions of groups carried out by  the exact sequence
 (\ref{centr}), together with its geometric interpretation, carries over to an arbitrary topos and 
therefore 
remains valid (except of the surjectivity of the right-hand arrow) 
when central extensions of topological groups, or of algebraic groups, are considered. 
The hypercohomology group $H^3(B, A)$, on the other hand, may no
 longer be described  by a short exact sequence
 (\ref{uct3}), since there now exists a new  non-trivial initial term in the spectral
sequence (\ref{hypcoh2}), provided by the group $\mathrm{Ext}^2(B, A)$. The  latter group was
 given a geometric interpretation
in \cite{pd:sga4} as the group of equivalence classes
of  strict Picard stacks. While this group vanishes in the category case, as we observed in 
\S \ref{sub:uct}, this is no longer true in  the general stack context.

\medskip
The alternating biextension point of view for analyzing $gr-$categories carries over
 very satisfactorily to the $gr-$stack context.  As we have already observed, 
 the 
 objects $X_x$, and  arrows (\ref{cxy})
 on which the definition of the 3-cocycle 
$f(x,y,z)$ depends no longer exist, but they do exist locally, so that the 3-cocycle $f(x,y,z)$ 
is locally defined. Alternate choices for these objects and arrows
yield cohomologous cocycles. Since the induced
map $\varphi_{\cc}$ (\ref{pic:wedgetrilin})  depends only on the cohomology class of $f$,  its
 local
representatives glue together to a globally defined  arrow form $\Lambda^3B$ to $A$.
Similarly, the weak biextension $E$ associated
to a $gr-$stack $\cc$ may be locally defined just as in (\ref{biexdef}), once representative
objects $X_x$ of $\cc$ have been chosen, and the local biextensions obtained in this manner glue
 to a weak biextension $E$ defined on all of $B \times B$.
 Its underlying torsor, however, is in general no longer endowed with a globally defined 
section $d$ (\ref{cxy1}),
 so that
the biextension $E$ may no longer be readily described in terms of cocycles. Propositions
\ref{pr:trilin} and \ref{pr:altern} remain valid in the stack context, and assert
 that the  weak biextension $E$ is a genuine (alternating)  biextension
of $B \times B$ by $A$ whenever the invariant $\varphi_{\cc}$ vanishes. This biextension may be
analyzed by the  methods of \cite{lb:alt}. As observed earlier, this yields a family of
induced   quadratic maps $\psi_n: {}_nB \la A $ whose vanishing  implies  by the universal
coefficient exact sequence (3.19) of \cite{lb:alt} that the biextension
$E$ descends
to an ordinary extension of $\Lambda^2B$ by $A$.

\medskip
Suppose now that the alternating biextension $E$  is 
trivial. We may then choose in a compatible manner,
 for each pair of sections $x,y$ of $B$, an
arrow 
\be
\label{ar:symO}
s(x,y):X_{y}X_{x} \la  X_{x}X_{y}
\ee
 in $\cc$. This actually
determines, for an arbitrary pair of objects $X,Y$ of $\cc$ (with
associated sections $x,y$ in $B \times B$), a symmetry  arrow  by the rule

\begin{equation}
\label{ar:sym}
\xymatrix{
{YX} \ar[r]^{c_y c_x}  & {X_y X_x} \ar[r]^{s(x,y)} & {X_{x}X_{y}}
\ar[r]^{c_{x}^{-1} c_{y}^{-1}} & {XY}
}
\end{equation}
and this is independent of the (local) choice of arrows  $c_x: X \la X_x$ and
$c_y: Y \la X_y$ in $\cc$. The compatibility of the section $s$ of $E$ with the partial
multiplication laws (\ref{eq: mult1}) and (\ref{eq: mult2}) implies, as we have already observed,
that
the symmetry arrows  (\ref{ar:symO}) (and therefore more generally the
corresponding symmetry arrows (\ref{ar:sym})) satisfy both hexagon
conditions, so that  $\cc$ is a braided stack.
Finally, compatibility of $s$ with the section $t$ (definition \ref{defalt})
of $E$ asserts that for $X = Y$ the composite map (\ref{ar:sym}) is simply the
identity map. This forces the braided category $\cc$ to be Picard strict, and so provides a direct 
geometric interpretation
of the degree 3 portion of the universal coefficient spectral sequence (\ref{hypcoh2}). 
We have therefore obtained in the stack context the following analog  of corollary \ref{cor:alt}.

\begin{corollary}
\label{cor:alt2}
  Let $\cc$ be a $gr-$stack  of $T$ satisfying the conditions of
proposition \ref{pr:altern}, and whose associated biextension $E$ is  
trivial (as an alternating biextension). Then $\cc$ is the underlying $gr-$stack
of a strict Picard stack.
\end{corollary}

Vanishing theorems for certain groups
$\mathrm{Ext}^{2}(B \, , \, A)$ are proved in  \cite{lb:invent}.
 When this vanishing takes place, the corresponding
strict Picard stack are trivial, so that the corollary \ref{cor:alt2}
 takes on the form of a vanishing 
theorem for $gr-$stacks, along the same lines as in corollary \ref{cor:alt}.

\section{Higher multiextensions}

\setcounter{equation}{0}%

\hspace{.7cm}The  cohomology groups $H^{n+1}(B, A)$ with $n > 2$ are the $k$-invariants of two-stage
 Postnikov
systems with homotopy groups $B$ and $A$  concentrated in degrees 1 and $n$. From a 
categorical point of view, such a cohomology class may  therefore be represented by a $(n-1
)-$category 
(or rather $(n-1)-$groupoid) $\cc$,  endowed with a multiplication $\cc \times \cc \la \cc$
 satisfying the requisite higher  associativity axiom. It is required that the group of isomorphism
classes of  objects of $\cc$ be isomorphic to $B$ and that the intermediate homotopy groups
$\pi_i(\cc)$ vanish ($0 < i < n-1$). The group $\pi_{n-1}(\cc)$, which is simply the group of self
$(n-1)-$maps of the identity $(n-2)-$arrow, is required to be isomorphic to $A$. Finally, the
requirement that
$A$ is a trivial
 $B$-module can be translated into a weak commutativity condition,
 analogous to hypothesis \ref{hyp1}.

\medskip
In the following discussion we will mainly be concerned with the  case
$n=3$, where the definition of a monoidal $2$-groupoid does not offer any difficulty. We note 
in passing that the analog  in this 2-categorical
  context  of a crossed module, which occurs
when the associativity isomorphism in $\cc$ is strict,  has been
worked out by D. Conduch\'e in \cite{mc.gen} definition 2.2 under the name of 2-crossed modules.
 The  requisite conditions on the homotopy groups now translates to the requirement that such a
2-crossed module $L \la M \la N$ lives in an exact sequence of groups
\[ 0 \la A \la L \la M \la N \la B \la 0 \]
A direct proof of the classification of such length 3 
extensions by
elements  of the group $H^4(B, A)$ is given  in \cite{mc.gen} theorem 4.7.
 Such a discussion can also  be carried
out from a categorical point of view by extending by one more step to a representation
of  pentagonal $2$-arrows the geometric construction of the 3-cocycle discussed in \S
\ref{sub:uct}.

\medskip
Let us now now examine the effect on cohomology  of the filtration on the chains on $K(B,1)$
by powers of the augmentation ideal
mentioned in the introduction. Recall that the terms
 which  occur in the analysis of  $H^{n+1}(B, A)$
are  the groups $\mathrm{Ext}^p(L\Lambda^qB, A)$, with $p+q = n+1$.
 In \mbox{particular,} the filtration on the group $H^4(B, A)$ (which
 classifies  $2$-categories  of the type
envisaged above)   yields successive geometric objects which respectively  live in  the groups
 $\mathrm{Hom}(\Lambda^4 B, A)$, $\mathrm{Ext}^1 (L\Lambda^3 B , A)$,
  $\mathrm{Ext}^2 (L\Lambda^2 B  , A)$
and  $\mathrm{Ext}^3 (B  , A)$. A prerequisite to  a geometric discussion of this filtration of $H^4$
is 
 the interpretation in geometric terms of the \noindent groups \(\mathrm{Ext}^p(\lotimes{}^{\!\!q}B , A) \)
for varying integers $p$ and $q$. We will call the objects  whose isomorphism classes 
are classified by 
these groups $(p,q)-$extensions (or $(p,q)-$multi-extensions) of $B$ by $A$.
 A $(p,1)$-extension is 
simply a $p$-fold extension by $A$ of the abelian group $B$, and is therefore geometrically
described by classes of strict Picard ($p-1$)-categories with invariants $B$ and $A$.
 Similarly, a $(1, 2)$-extension is an
(ordinary) biextension of $B \times B$ by $A$ and more generally a $(1, q)$-extension is,
 in the terminology of \cite{ag:sga7} VII 2.10.2, a
 $q$-extension of the $q$ groups $B_1, \ldots , B_q$ by $A$, with $B_1 = \cdots = B_q = B$.
  These  interpretations  of $(p,1)$ and $(1,q)$-extensions may be combined as follows. Choosing
 in the manner explained
in \cite{ag:sga7} VII as
representative for the object $L\otimes^q B $ of the derived category   the $q$-fold tensor
product of a canonical free resolution of $B$,  it is apparent that  a $(p,q)-$extension for a
general  pair of integers $p$ and $q$ may be thought of as an abelian $(p-1)-$gerbe $\cc$ on $B^q$
 \cite {lb:2-gerbes},
together with a family of $q$ partial group laws
 $\stackrel{1}{+}, \cdots , \stackrel{q}{+}$ on $\cc$
living above the  $q$ composition laws  on $B^q$ determined by the group laws on each of the 
$p$ factors. Each of these partial group laws is required to  satisfy
 the requisite higher associativity
and commutativity conditions, together with  higher compatibily conditions between them. These
higher conditions may be worked out by considering the cells and their boundaries in the chosen
representative of $L\!\otimes^q \!B $. We simply spell this out in the case of ($2,2$)-extensions,
the  only essentially new case required for an understanding of   $H^4(B, A)$. As we have just
asserted, this is an abelian $A$-gerbe $\cc$ above $B \times B$, together with a pair of partial
group laws
 $\stackrel{1}{+}$
and $\stackrel{2}{+}$. The partial commutativity and associativity conditions assert that,
 for each section $x: S \la B$ of $B$, the groups laws $\stackrel{2}{+}$
and $\stackrel{1}{+}$ respectively endow the pullbacks $(x \times 1)^{\ast}\cc$ and   
$(1 \times x)^{\ast}\cc$  of $\cc$ above the $S$-groups $S \times B$ and $B \times S$ with the
structure of a  strict Picard stacks \cite{lb:2-gerbes}. The compatibility conditions between
the
two group laws are described as follows. We may choose functorial isomorphisms 
\be
\label{comp12}
c_{x_1, x_2 \, ; \, x_3, x_4}: (X \oa Y) \ob (Z \oa W) \la (X \ob Z) \oa (Y \ob W)
 \ee where the projections $\pi$ of the  four objects
$X,Y,Z,W$  to the group of isomorphism classes of
objects satisfy 
\[
\begin{array}{cccc}
\pi(X) = (x_1, x_3) \: \:  & \pi(Y) = (x_2, x_3) \: \:  & \pi(Z) = (x_1, x_4) \: \: & \pi(W) = (x_2,
x_4)
\end{array}
\]
for  sections  $x_i $ of $B$, so that the source and target of (\ref{comp12}) are well defined.
These isomorphisms $c$ are required to be compatible with the associativity and commutativity 
isomorphisms
 for $\oa$ and $\ob$. The compatibility of $c$ with the commutativity isomorphisms
asserts that
the following diagram, in which the horizontal arrows are
determined by the commutativity axiom for $\ob$,   commutes for all allowable object. So must the
corresponding one in which the role of $\oa$ and $\ob$ have been exchanged.
\begin{equation}
\label{2commut}
\xymatrix{
{(X \oa Y) \ob (Z \oa W)} \ar[r] \ar[d]_>>>>>c &  {(Z \oa W) \ob (X \oa Y)} \ar[d]^>>>>>c\\
{(X \ob Z) \oa (Y \ob W)} \ar[r] & {(Z \ob X) \oa (W \ob Y)}
}
\end{equation}
 Similarly, the compatibility of the maps (\ref{comp12}) with the associativity isomorphisms for the
two partial group laws  is described by the commutativity of the following diagram in which the
 vertical arrows are  maps (\ref{comp12}) and the horizontal diagrams are
associativity isomorphisms for $\ob$, and by the corresponding one in which the role of $\ob$ and
$\oa$ is interchanged.

\begin{equation}
\label{compass2}
\xymatrix{
{(X \oa Y) \ob ((Z \oa W) \ob (S \oa T))}\ar[r]\ar[d]
&{((X \oa Y) \ob (Z \oa W)) \ob (S \oa T)}
\ar[d]\\
{(X \oa Y) \ob ((Z \ob S) \oa (W \ob T))} \ar[d] & {((X \ob Z) \oa (Y \ob W)) \ob (S \oa T)} \ar[d]
\\{(X \ob (Z \ob S)) \oa (Y \ob(W \ob T))} \ar[r] &
{((X \ob Z) \ob S) \oa ((Y \ob W) \ob T)}}
\end{equation}

\medskip
\noindent When   associativity and
compatibility constraints  are  given for the laws $\oa$
 and $\ob$, but no  commutativity constraint, 
(and {\it a fortiori} no commutative diagram (\ref{2commut}) is introduced), we will say that
$\cc$ is a weak (2,2)-extension of $B \times B$ by $A$.  Observe that, despite their somewhat
abstract aspect, (2,2)-extensions are not hard to classify. Indeed, the group $\mathrm{Ext}^2 (B
\lotimes B, A)$ of  isomorphism classes of such (2,2)-extensions may be analyzed {\it via} the
adjunction spectral sequence, whose low-degree terms are described in \cite{ag:sga7} VIII (1.1.4).
In that context, when $B$ is an abelian variety over an algebraically closed field,
and $A$ is the multiplicative group $G_m$, the vanishing of both $\mathrm{Hom}(B,A)$ and
$\mathrm{Ext}^2(B, A)$ (see \cite{lb:invent}) ensures that such (2,2)-extensions are classified
up to equivalence by the group $\mathrm{Ext}^1 (B,B^t)$ of extensions of $B$ by the dual abelian
variety
$B^t$.

\medskip
In order to
describe alternating (2,2)-extensions, we  first introduce the concept of an anti-symmetric
(2,2)-extension. In the situation just examined, these have a very concrete interpretation,
since they are described by  extensions
\be
\label{Bt}
0 \la B^t \la E \la B \la 0
\ee
of
$B$ by
$B^t$ which are opposite (for the Baer sum) to the  extension
\[ 0 \la B^t \la E^t \la B \la 0\]
obtained by applying the contravariant duality functor $( \: )^t$ to the exact sequence 
of abelian varieties (\ref{Bt}).
 Returning to the general situation, these anti-symmetric (2,2)-extensions may be described as
follows. Let
$s$  once more denote the map which permutes the factors of
$B^2$.
\begin{definition}
\label{2asdef}
A (2,2)-extension is anti-symmetric if it is endowed with
 a morphism of (2,2)-extensions
\be
\label{pi2}
 \pi: \cc^{-1} \la s^*\cc
\ee
analogous to isomorphism (\ref{pialt}),  whose source is the $A$-gerbe $\cc^{-1} = \mathrm{Hom}_A
(\cc , \mathrm{Tors} A)$ of morphisms\footnote{A morphism $\Phi: \cc \la \mathcal{D}$ of $A$-gerbes
is a morphism of gerbes for which the induced maps $Aut(X) \la Aut (\Phi X)$ are identified with the
identity by the
$A$-gerbe structure on the source and target gerbe. In particular,
such a morphism $\Phi$ induces (in the terminology of Giraud \cite{jg:cohomologie}) the identity map
on liens.} of
$A$-gerbes from
$\cc$ to the trivial
$A-$gerbe.  We further  require that  the composite morphism
\be
\label{pi3} s^*\pi_{\cc}  \circ \pi_{{\cc}^{-1}}: \: \cc \simeq (\cc^{-1})^{-1} \la s^*\cc^{-1} \la
s^*s^*\cc
\simeq
\cc
\ee
 be equivalent to the identity functor, by an equivalence which is unchanged when the factors of
$B^2$ are permuted. 
\end{definition}

 For any stack $\cc$, let us
denote  by $\cc^{0}$ the opposite stack of $\cc$, whose fibers are the  categories  opposite to the
fibers of $\cc$. If
$\cc$ is a gerbe, then it is immediate that $\cc^{0}$  also is  one. Suppose further that $\cc$ is an
$A$-gerbe for some group $A$, so that we are given, for objects $X$ of $\cc$, a family of isomorphism
$\lambda_X: A \la Aut_{\cc}(X)$. Then these maps $\lambda$ may also be viewed as isomorphisms
$\lambda_X: A^{0}
\la Aut_{\cc^{0}}(X)$ between the opposite groups, so that  they define on $\cc^{0}$ a natural 
$A^{0}-$gerbe structure.  We believe that the following description
of the inverse  $\cc^{-1} $ of an abelian $A$-gerbe $\cc$ may be of independent interest. It is a
local statement,  and may therefore  be  verified by supposing that $\cc$ is the trivial $A-$gerbe,
in which case it is immediate.

\begin{lemma}
Let $A$ be an abelian group of $T$, and $\cc$ an abelian $A$-gerbe. The Yoneda morphism
\[
\begin{array}{ccc}
\cc^{0} & \la & \mathrm{Hom}_A (\cc ,\mathrm{Tors} A) \\
 X &\longmapsto & h^{X}
\end{array}
\]
is an isomorphism of $A$-gerbes. 
\end{lemma}

The morphism (\ref{pi2}) which defines an anti-symmetry structure on the $(2,2)-$extension $\cc$ 
may therefore be described by a morphism of $A-$gerbes
\be
\label{contrpi} \pi_{\cc}: \cc^{0} \la s^*\cc 
\ee 
in other words  by a ``contravariant morphism of
$A$-gerbes'' from $\cc$ to its pullback $s^* \cc$ which is compatible, in the obvious sense,  with
each of the two given partial group laws on $\cc$. We further  require that the additional
anti-symmetry conditions on
$\pi_{\cc}$ analogous to those of definition \ref{2asdef} be satisfied. The pullback $\Delta\cc$
of $\cc$ along the diagonal is then canonically endowed with the structure of a strict Picard
category with invariants
$B$ an $A$, and in fact one whose pullback by the ``multiplication by 2'' map (\ref{mult2}) is   a
trivial  strict Picard category.  The full definition of an alternating
$(2,2)$-extension may now be given.

\begin{definition}
\label{22alt}
A $(2,2)-$extension $\cc$ of $B \times B$ by $A$ is alternating if it anti-symmetric, and if the
induced strict Picard category 
$\Delta
\cc$ is trivial, by a trivialization whose pullback   by the ``multiplication by 2'' map
(\ref{mult2}) is the canonical one determined by the anti-symmetry structure\footnote{It is
equivalent to require that the square of the trivialization is that on  $\cc^2$ determined by the
anti-symmetry structure.}.
\end{definition}

The other higher element of structure whose definition  will be  required is the concept of an
alternating (1,3)-extension, an object classified up to equivalence by elements of the group
$\mathrm{Ext}^1(L\Lambda^3B, A)$. We begin with an ordinary
(1,3)-extension (in other words a triextension) of $B^3$ by $A$. Recall that this is an 
$A$-torsor $E$ on $B^3$, together with three partial multiplication laws
 \[
\begin{array}{cccc}
\stackrel{1}{+}:   & E_{x,y,z} \: E_{x',y,z}  & \la & E_{x+x', y,z} \\ \stackrel{2}{+}: & E_{x,y,z} \:  E_{x,y',z} & \la & 
E_{x,y+y',z} \\
 \oc: & E_{x,y,z} \:  E_{x,y,z'} &  \la  & E_{x,y,z+z'}
\end{array}
\] 
each
of which is commutative and associative, and any two of which are compatible with each other. Just
as  alternating biextensions could be understood by considering the  complex (\ref{kosl2}),
information concerning an alternating structure on $E$ may be inferred from the  complex
\be
\label{kosl3}
 0 \la \Gamma_3B \la \Gamma_2 B \otimes B \la B \otimes \Lambda ^2 B \la \Lambda^3 B \la
0 
\ee
which is simply a Koszul complex \cite{li: cotg1} chapter 1 (4.3.1.3). By examining the right-hand
arrow in this complex, we see that an alternating triextension $E$ with general fibre $E_{x,y,z}\,$,
when viewed for a fixed
$x \in B$  as a biextension in $y,z$ of $B^2$ by $A$, must  be alternating in the
variables $y,z$. Its restriction $E_{x,z,z}$ above the diagonal $\Delta_{23}$ is therefore
provided with a section
$t^2_{x,z}
\in E_{x,z,z}$. It is required that this section be compatible, as in definition \ref{defalt},
with the symmetry and cube structures in $z$ on $E_{x,z,z}$  determined by the second and third
group laws on 
$E$. This section $t^2_{x,z}$ must also be linear in $x$, in the sense that the equation
 \be
\label{compx}
 t^2_{x,z} \stackrel{1}{+} t^2_{x',z} = t^2_{x+x',z}
\ee
must be satisfied in $E_{x+x',z}$. The middle map in  the sequence
(\ref{kosl3}) may be interpreted as determining a constraint on the triextension $E$ which is
 to be quadratic in $x$ and linear in $z$. Such a constraint is given for a fixed $z$  by
an alternating structure on $E_{x,y,z}\,$, viewed  as a biextension in the variables
$x,y$, in other words by a section  
$t^1_{x,z} \in E_{x,x,z}$ compatible with the $\Sigma$ -structure in $x$ and satisfying  a
relation analogous to (\ref{compx}) in the variable $z$. 

\medskip
We must now express the
compatiblity conditions between these sections $t^1$ and $t^2$ which follow from the left-hand arrow
in  (\ref{kosl3}). Since the group $\Gamma_3B$ has two separate sorts of generators, those of the
form $\gamma_3(b)$ and those of type $\gamma_2(b) \, b'$ for elements $b,b' \in B$, two distinct sorts
of compatibilities will have to be verified between these sections. These compatibilities will
however be related to each other by the conditions corresponding to the identities 
\begin{eqnarray*}
3\gamma_3(b) & = & \gamma_2(b)b \\ 
\gamma_3(b+b') - \gamma_3(b) - \gamma_3(b') & =& \gamma_2(b)\, b' + \gamma_2(b') \, b
\end{eqnarray*} in the group
$\Gamma_3B$. The first of these compatibilites, which is cubical in $b$
and therefore corresponds to the generators  $\gamma_3(b)$ of  $\Gamma_3B$,
asserts that for all
$x \in B$ the equation
\begin{eqnarray}
\label{tri1}
t^1_{x,x}& = & t^2_{x,x}
\end{eqnarray}
is satisfied in $E_{x,x,x}$. In order to state the second compatibility condition between $t^1$ and
$t^2$ as pleasantly as possible, we need the following lemma.

\begin{lemma}
\label{pr:alttriext}
Let $E_{x,y,z}$ be a triextension of $B^3$ by $A$, endowed  as above with an alternating
structure in $x$ linear in $z$ defined by a section
$t^1_{x,z} \in E_{x,x,z}\,$, and with an
alternating structure in $z$ linear in $x$ defined by a section
$t^2_{x,z} \in E_{x,z,z}$. The restriction $E_{x,y,x}$
of $E$ above the diagonal $\Delta_{13}$ is canonically endowed with an alternating structure
$t^3_{x,y}
\in E_{x,y,x}$ in $x$,  which is linear in $y$.
\end{lemma}
{\bf Proof:} Since $E$ is alternating with respect to $x$, it is {\it a fortiori} anti-symmetric. This
amounts to the asssertion that the section 
\(
t^1_{x+y,z} \: (t^1_{x,z})^{-1} \: (t^1_{y,z})^{-1} \in E_{x+y, x+y,z}\:  E^{-1}_{x,x,z} \:
E^{-1}_{y,y,z}
\)
defines, {\it via} a canonical
isomorphism determined by the biextension structure of $E$, a trivialization
$s^1_{x,y,z}$ of the symmetric biextension 
$
F^1_{x,y,z}\:~=~\: E_{x,y,z} \: E_{y,x,z}
$. Since $t^1$ defines a
$\Sigma$-structure in $x$ on  $E_{x,x,z}$, and is therefore  quadratic in the variable
$x$, the equation
$t^1_{2x,z}~=~(t^1_{x,z})^4$ is satisfied up to canonical isomorphism, from which the equation
\begin{equation}
\label{s1t1}
 s^1_{x,x,z} = (t^1_{x,z})^2 
\end{equation}
in $E_{x,x,z}$ follows immediately. The equation
\begin{equation}
\label{s2t2}
s^2_{x,y,y} = (t^2_{x,y})^2
\end{equation}
(where $s^2_{x,y,z} = t^2_{x,y+z} \: (t^2_{x,y})^{-1} \: (t^2_{x,z})^{-1}$) is proved in the same way,
as a consequence of the quadraticity in
$z$ of
$t^2_{x,z}$. Let us now  set 
\begin{eqnarray}
\label{t31} 
t^3_{x,y} & = & s^1_{x,y,x}(t^2_{y,x})^{-1}
\end{eqnarray}
This is, as required, an element of $(E_{x,y,x} \:  E_{y,x,x})(E_{y,x,x})^{-1} \simeq E_{x,y,x}$. It
is readily verified that the section  $t^3_{x,y}$ of $E_{x,y,x}$ defined in this manner satisfies the
 requisite quadraticity condition in $x$  and linearity condition in $y$, so that the lemma is proved.

\medskip
Another section $\tilde{t}\,^3_{x,y}$ of $E_{x,y,x}$ with the same properties  as $t^3_{x,y}$ could
have
 been defined in terms of the alternating structure $t^1_{x,z}$ and the antisymmetry
section
$s^2_{x,y,z} \in E_{x,y,z} \: E_{x,z,y}$ determined by $t^2_{x,y}$ by setting 
\begin{eqnarray}
\label{t32} 
\tilde{t}\,^3_{x,y} & = & s^2_{x,y,x}(t^1_{x,y})^{-1}
\end{eqnarray}
The second compatibility which the sections $t^1$ and $t^2$ must satisfy is the
requirement that
\begin{eqnarray}
\label{t33} 
t^3_{x,y} & = & \tilde{t}\,^3_{x,y}
\end{eqnarray}
This may, of course, also be written as the condition
\begin{eqnarray}
\label{t4}
\frac{t^1_{y,x}}{ t^2_{x,y}} & = & \frac{s^2_{y,x,y}}{s^1_{y,x,y}}
\end{eqnarray}

Our definition of an alternating triextension is now complete. It can be summarized as follows, with
the corresponding notion of a trivialization spelled out.

\begin {definition}
\label{def:triex}
A triextension $E$ of $B^3$ by $A$ is alternating  if it is endowed with sections
$t^1_{x,z}
\in E_{x,x,z} $ and $ 
t^2_{x,z} \in E_{x,z,z}$ such that $t^1$  defines a partial alternating structure on $E$ with
respect to $x$ linear in $z$ and $t^2$ defines a partial alternating structure on $E$ with
respect to $z$ linear in $x$. The sections $t^1$ and $t^2$
must  also satisfy the compatibility  conditions (\ref{tri1}) and (\ref{t4}) with $s^1$ and $s^2$
defined as in the proof of lemma \ref{pr:alttriext}. A trivialization of
$E$ as an alternating triextension is determined by a section $\sigma_{x,y,z}$ of $E_{x,y,z}$ which
trivializes $E$ as a triextension (in other words compatibly with each of the three partial group
laws), and such that 
\begin{eqnarray*} 
\sigma_{x,x,z} & = & t^1_{x,z} \\
\sigma _{x,z,z} & = & t^2_{x,z}
\end{eqnarray*}
\end{definition}

The section $t^3_{x,y}$ of such an alternating triextension defined by formula (\ref{t31}) determines
as above  a partial antisymmetry structure $s^3_{x,y,z} \in E_{x,y,z}\:E_{z,y,x}$, by setting, up to a
canonical isomorphism
\begin{equation}
s^3_{x,y,z}   =  \frac{t^3_{x+z, y}}{t^3_{x,y}t^3_{z,y}} 
\end{equation}
The formula 
\be
\label{s123}
s^3_{x,y,z} = \frac{s^1_{x,y,z}s^1_{z,y,x}}{s^2_{y,x,z}}
\ee
now follows from (\ref{t32}) and the definition of $s^2_{x,y,z}$.
The equations 
\begin{eqnarray}
\label{t123}
t^1_{x,z} & = & s^3_{x,x,z} (t^2_{z,x})^{-1} \nonumber \\
t^2_{x,z} & = & s^3_{x,z,z} (t^1_{z,x})^{-1}
\end{eqnarray}
are  consequences of (\ref{s1t1}), (\ref{s2t2}) and (\ref{t4}). The first of these equations 
shows  that this new method for constructing
$t^1$ out of
$t^2$ and the anti-symmetry condition $s^3$ derived from $t^3$ yields the same result as the
method (\ref{t32}) for constructing $t^1$ out of $t^3$ and the anti-symmetry  condition $s^2$ derived
from
$t^2$. Similarly the second formula shows that the new method (\ref{t123}) for constructing $t^2$
out of $t^1$ and the anti-symmetry condition derived from $t^3$ yields the same result as the method
(\ref{t31}) for constructing it out of $t^3$ and the anti-symmetry condition $s^1$ derived from
$t^1$. Allowing ourselves a certain amount of redundancy, we may therefore give another  description
of alternating triextension which is entirely symmetric in the variables $x,y,z$, as befits an
object associated to $\Lambda^3B$.

\begin{proposition}
A triextension $E$ of $B^3$ by $A$ is alternating if and only if it is endowed with sections
$t^1_{x,z}
\in E_{x,x,z} \, , \, 
t^2_{x,z} \in E_{x,z,z}$ and
$t^3_{x,y} \in E_{x,y,x}$ each of which defines a partial alternating structure with
respect to the repeated variable which is linear with respect to the other variable and which 
 satisfy the following compatibility conditions
\begin{enumerate}
\item For each $i$, the two possible methods  described above for constructing a section $t^i$
in terms of the two other sections $t^j$ and $t^k$ yield the same result.
\item For
every $x \in B $, the equation $t^1_{x,x} \: = \: t^2_{x,x} \: = \: t^3_{x,x}$ is satisfied in
$E_{x,x,x}$.
\end{enumerate}
\end{proposition}

\begin{remark} {\rm This description of alternating triextensions may  be obtained  in a
somewhat more symmetric manner by making use of the
derived version of the sequence
\[
0 \la \Gamma_3B \la (\Gamma_2 B \otimes B) \oplus (B \otimes \Gamma_2 B) \la B^{\, \otimes 3}
\la
\Lambda^3 B
\la 0 
\]
 of \cite{akin}, instead of the Koszul sequence (\ref{kosl3}). This sequence also makes it
immediately clear that an alternating triextension whose underlying triextension is trivial, may be
described in terms of pairs of compatible maps $f,g: B \times B \la A$, with $f$ quadratic in the
first variable and linear in the second one ({\it resp.} g linear in the first variable  and 
quadratic in the second one). Note that in an algebro-geometric setting, this often implies that
such alternating triextensions are trivial. for example, when $B$ is an abelian variety over an
algebraically closed field, it is easily verified (see \cite{ag:sga7} VII 2.10.2) that any
triextension of $B$ by the multiplicative group $G_m$ is trivial. The assertion is now immediate,
since the only maps from $B \times B$ to $G_m$ are the constant ones. }
\end{remark}

\section{Picard structures on monoidal 2-categories}

\setcounter{equation}{0}%

\hspace{.7cm}Let $\cc$  be a monoidal group-like 2-stack in groupoids,  as defined under the
name of 2-{\it gr}-stack  in \cite{lb:2-gerbes} definition 8.4.  The first invariant  of $\cc$
 is the sheaf of groups  $\pi_0(\cc)$ associated to the presheaf of 
isomorphism classes of objects of
$\cc$.  In the category case which we will mostly consider, we will  say that $\cc$ is a
monoidal group-like 2-groupoid.  
The group $\pi_0(\cc)$ is then simply the group of isomorphism classes of objects of
$\cc$. Our first assumption will be, as in the monoidal 1-category case,
that this group  is abelian.  Since the  monoidal category
$\mathrm{Aut}_{\cc}(I)$ of self-arrows of the unit object $I$ of $\cc$ is automatically braided, the
two other homotopy groups of
$\cc$,  which may  be defined by 
\be
\pi_i(\cc) = \pi_{i-1}(\mathrm{Aut}_{\cc}(I))
\ee
for $i= 1,2$, are both abelian groups. It was explained in \cite{lb:2-gerbes}  that such
monoidal 2-groupoids with $
\pi_0(\cc) 
\simeq B$ and
$\mathrm{Aut}_{\cc}(I)$ equivalent to a given braided category $\mathcal{A}$ are classified by an
appropriately defined cohomology group $H^3(B, \mathcal{A})$. The group $B$ acts by conjugation on
the category $\ac$, and we will assume that this action is equivalent to the trivial one\footnote{in
other words that  the tensor functor $\phi: \underline{B} \la \mathrm{Aut}_{\cc}(I)$
with source the discrete category defined by $B$ is equivalent to the trivial functor.}. Finally, we
will assume in the sequel for simplicity that the abelian group
$\pi_1(\cc) $ is trivial, so that
$\ac$ is the category with a unique object whose arrows form an abelian group $A$. The cohomology
group  $H^3(B, \mathcal{A})$ then reduces to the standard cohomology group
$H^4(B, A)$ with values in the trivial $B$-module $A$. As we have said, the class of $\cc$ may be
viewed as the
$k$-invariant of the two stage Postnikov system 

\[
\xymatrix
{{K(A,3)}\ar[r]&{X}\ar[d]\\
& {K(B, 1)}}
\]  
defined by  the the  classifying space $X$ of
 the nerve of the monoidal 2-category $\cc$.  In more explicit terms, one associates to $\cc$ the
$A$-valued 4-cocycle
$f(x,y,z,w)$ obtained as follows. Choose, as in the case of monoidal categories, representative
objects
$X_x$ and arrows
$c_{x,y}: X_x X_y \la X_{xy}$ in
$\cc$. Since it is assumed here that $\pi_1(\cc) = 0$, we may also choose for every $x,y,z \in
B$ a 2-arrow 
\[  \eta_{x,y,z}:1_{X_{x,y,z}}  \Longrightarrow f(x,y,z) \]
between the identity 1-arrow, and the 1-arrow defined as in \ref{diag:f}. The
pentagon 2-arrow associated to the four objects
$X_x,
\,  X_y,
\, X_z,
\, X_w$ then determines an element $f(x,y,z,w)$ in $\mathrm{Aut}_{\cc}(X_{xyzw})$, {\it i.e.} a
4-cochain $f:B^4 \la A$. Stasheff's $K_5$ relation \cite{sta} implies that $f$ is a 4-cocycle. Other
choices of objects $X_x$, 1-arrows $c_{x,y}$ and 2-arrows 
$\eta_{x,y,z}$  determine cohomologous cocycles so that the class of $f$ in $H^4(B,A)$ only depends
on the equivalence class of $\cc$.

\medskip
As in our study of monoidal categories, we may analyze  the monoidal 2-category $\cc$ by
introducing, for each pair of elements $x, y \in B$, the commutator category
\be
\label{2mock}
 \mathcal{E}_{x,y} = Isom_{\cc} (X_y X_x, X_x X_y)
\ee
This is a groupoid, on which the $gr$-category of self-arrows of $X_x X_y$ acts  on the right
fully and  faithfully by
composition of arrows. This {\it}-category is equivalent to $\mathcal{A}$, and therefore,  since
$\pi_1(\cc) = \mathrm{0}$, to the groupoid 
$A[1]$ with a single object defined by the abelian group $A$. The categories $ \mathcal{E}_{x,y}$ 
assemble, for varying
$x,y$, to form an abelian $A$-gerbe $\mathcal{E}$  on $B \times B$ which is a first element of
structure associated to the monoidal 2-stack $\cc$. In fact this
gerbe is trivial, since one can choose a compatible family of   sections $s(x,y) \in 
\mathcal{E}_{x,y}$, for example as in (\ref{cxy1}) those obtained by composing the chosen maps
$X_x X_y \la X_{xy}$ with an inverse of the maps $X_y X_x \la X_{yx}$. There remains however some
interesting structure on $\mathcal{E}$ to be explored.
Indeed, the constructions (\ref{u+1v}) and (\ref{u+2w}) now define  partial group laws  $\stackrel{1}{+}$ and
$\stackrel{2}{+}$ on  the abelian $A$-gerbe $\mathcal{E}$ on
$B \times B$. Once more, no commutativity property  for the partial laws $\stackrel{1}{+}$ and
$\stackrel{2}{+}$ is asserted, so that $\mathcal{E}$ is in general a weak, rather than a genuine
(2,2)-extension.

\medskip
This analysis of $\cc$, and of its associated commutator $\mathcal{E}$ carries over from
monoidal 2-categories to monoidal 2-stacks, the only significant difference being that in that case
the underlying $A$-gerbe of
$\mathcal{E}$ is no longer trivial. The following higher analog of proposition
\ref{pr:Eweak} is therefore true.
\begin{proposition}
\label{weak22ext}
Let $\cc$ be a monoidal 2-stack with invariants $B$ and $A$, satisfying the previous hypotheses.
$\cc$ is classified up to equivalence by an element of the (hyper)-cohomology group $H^4(B, A)$ (for
A a trivial
$B$-module). The  constructions (\ref{u+1v}) and (\ref{u+2w}) define  on  the abelian $A$-gerbe
$\mathcal{E}$ (\ref{2mock}) on 
$B \times B$ a weak (2,2)-extension
 structure. 
\end{proposition}

In order to prove  that each of the two  group laws on $\cc$ is coherently associative, one could
simply examine the next higher versions of diagram (\ref{pic:assoc1}), in other words the weak
versions of the pair of  2-categorical diagrams which would, in the terminology of \cite{kv:2-cat}, be
denoted by
$(\bullet
\otimes (\bullet
\otimes \bullet \otimes
\bullet
\otimes
\bullet))$ and $((\bullet
\otimes \bullet \otimes \bullet \otimes \bullet) \otimes \bullet)$.  The next higher version of the
compatibility diagram (\ref {pic:compat1}) would then show that the  two group laws are compatible
with each other, in the sense made explicit for the map (\ref{comp12}) by the two diagrams
(\ref{compass2}). Such an argument would certainly be sufficient in order to prove  the proposition. 
However, if one wanted to fill in the details of a proof along these lines, one would be led to the
consideration of a family of commuting 2-categorical diagrams, which cannot be represented here in
an enlightening manner. We therefore prefer to give a proof of the proposition  in cocyclic, rather
than diagrammatic  terms, even though this method of proof {\it a priori} only  applies in the
monoidal 2-category case, rather than the full monoidal 2-stack situation. The method of proof which
we now propose will thus be analogous to the discussion in remark
\ref{cocyclbiext}, but at the next higher level.

\medskip

Starting from an $A$-valued 4-cocycle $f(x_1,x_2,x_3,x_4)$, we have  seen  that for a fixed $x \in
B$ the group law $\stackrel{1}{+}$ is obtained by inserting  2-arrows derived from $f$ into the pentagons
 by which the vertices of diagram (\ref{pic:assoc1}) were replaced when we passed from this diagram
to its non-strict version. The 3-cochain
$\psi_{x_4}: B^3
\la A$ which describes as in (\ref{diag:f}) the associativity morphism in the monoidal category
$\mathcal{E}_{(
\:, x_4)}$ on
$B$ for the group law $\stackrel{1}{+}$ is therefore defined by composition of these 2-arrows, in other words
(once the sign has been taken into account) by the formula
\be
\label{psiass}
\psi_{x_4} (x_1,x_2,x_3) = \prod_{\sigma(1) < \sigma(2) <\sigma(3)} f(x_{\sigma(1)} , x_{\sigma (2)},
x_{\sigma (3)} , x_{\sigma (4)})^{- \epsilon(\sigma)}
\ee
in which $\epsilon (\sigma)$ denotes the sign of the permutation $\sigma$. This is just the
product of the signed permutations of
$f(x,y,z,w)$ when the variable
$w$ is shuffled through
$(x,y,z)$. The group law $\stackrel{2}{+}$ on $\mathcal{E}_{(x_1, \:)}$ is similarly described  by 
\be
\label{phiass}
 \phi_{x_1} (x_2,x_3,x_4) = \prod_{\sigma(2) < \sigma(3) <\sigma(4)} f(x_{\sigma(1)} , x_{\sigma (2)},
x_{\sigma (3)} , x_{\sigma (4)})^{\epsilon(\sigma)}
 \ee
{\it i.e.} by the product of the signed shuffles in $f$ of  $x_1$ 
through
$(x_2,x_3,x_4)$. That each of the two group laws satisfies the pentagon condition is equivalent to the
assertion that the corresponding 3-cochain (\ref{psiass}),(\ref{phiass}) is a 3-cocycle, and this
follows readily from the 4-cocycle condition on $f$. In fact, it is unnecessary to perform this
computation explicitly, in view of the following observation. Consider Eilenberg-Mac Lane's iterated
bar-construction model $A(B,2)$ 
\cite{eml:hpin} \S 14 for the complex  of chains on the Eilenberg-Mac Lane
space $K(B,2)$. Since this is a chain complex, the square $\delta \circ \delta$ of the differential
$\delta$  is trivial when applied to any cell $c$. Applying this respectively to the cells  $[x_1
\mid_2 x_2,x_3,x_4,x_5]$ and
$[x_1,x_2,x_3,x_4
\mid_2 x_5]$, and passing from chains to
$A$-valued cochains on $K(B,2)$  yields the sought-after assertion. With this in mind, we relabel the
two  previous associativity maps by setting
\[
\begin{array}{ccc}
\phi (x_1,x_2,x_3\mid_2 x_4) & = & \psi_{x_4} (x_1,x_2,x_3)\\
\phi (x_1 \mid_2 x_2,x_3,x_4) & = &  \phi_{x_1} (x_2,x_3,x_4)  
\end{array}
\]
even though the first of these definitions is only consistent with \cite{eml:hpin} up to a sign.
Similarly, the compatibility isomorphism (\ref{comp12}) between the two group laws on
$\mathcal{E}$  is described by the cochain 
\[
\phi (x_1,x_2 \mid_2 x_3,x_4) = \prod_{\sigma(1) < \sigma(2)\: ; \: \sigma(3) <\sigma(4)}
f(x_{\sigma(1)} , x_{\sigma (2)}, x_{\sigma (3)} , x_{\sigma (4)})^{- \epsilon(\sigma)}
\]
obtained by shuffling $(x_1,x_2)$  through $(x_3,x_4)$. The vanishing of the image under
 $\delta \circ \delta$ of the cells $[x_1,x_2,x_3 \mid_2 x_4,x_5]$  and $[x_1,x_2 \mid_2
x_3,x_4,x_5]$ (or a direct computation)  imply that the higher compatibility conditions
(\ref{compass2}) are satisfied in
$\mathcal{E}$, so that  proposition \ref{weak22ext} is proved.

\medskip

In order to understand under which conditions the weak monoidal commutator (2,2)-extension obtained
from proposition \ref{weak22ext} is a genuine (2,2)-extension (in other words one whose partial
group laws are strictly commutative), we need only apply to the monoidal categories 
$\mathcal{E}_{(
\:, x_4)}$ and $\mathcal{E}_{(x_1, \:)}$  the theory developed in sections
\ref{sec:comweak}-\ref{par: altcom}. By (\ref{biext:1})-(\ref{biext:2}), the weak biextension
associated to the monoidal category $\mathcal{E}_{(
\:, x_4)}$ is described,  for a fixed $x_4 \in B$, by the  cochains 
\begin{eqnarray}
g(x_1,x_2;x_3 \mid_2 x_4) & = & \prod_{\sigma (1) < \sigma (2)} \phi (x_{\sigma (1)},
x_{\sigma (2)}, x_{\sigma (3)} \mid_2 x_4)^{\epsilon (\sigma)} \\
& = & \prod_{\sigma(1) < \sigma(2)}  f(x_{\sigma(1)} , x_{\sigma (2)},
x_{\sigma (3)} , x_{\sigma (4)})^{- \epsilon(\sigma)} \nonumber
\end{eqnarray}
and 
\begin{eqnarray}
h(x_1;x_2,x_3 \mid_2 x_4) & = & \prod_{\sigma (2) < \sigma (3)} \phi (x_{\sigma (1)},
x_{\sigma (2)}, x_{\sigma (3)} \mid_2 x_4)^{- \epsilon (\sigma)} \\ 
& = & \prod_{\sigma(2) <\sigma(3)}
 f(x_{\sigma(1)} , x_{\sigma (2)}, x_{\sigma (3)} , x_{\sigma(4)})^{\epsilon (\sigma)} 
\nonumber
\end{eqnarray} 
Similarly, the weak biextension associated  to the monoidal category $\mathcal{E}_{(x_1,
\:)}$ is described, for a fixed $x_1 \in B$, by the pair 
\begin{eqnarray*}
\gamma(x_1 \mid_2 x_2,x_3;x_4) & = & \prod_{\sigma (2) < \sigma (3)} \phi(x_1 \mid_2
x_{\sigma (2)}, x_{\sigma (3)}, x_{\sigma (4)})^{\epsilon (\sigma)} \\
\eta(x_1 \mid_2 x_2; x_3, x_4) & = & \prod_{\sigma (3) < \sigma (4)}
\phi (x_1 \mid_2 x_{\sigma (2)}, x_{\sigma (3)},x_{\sigma (4)})^{- \epsilon (\sigma)}
\end{eqnarray*}
so that 
\[\gamma(x_1 \mid_2 x_2, x_3; x_4) = h(x_1;x_2,x_3 \mid_2 x_4)
\]
and 
\begin{eqnarray}
\eta (x_1 \mid_2 x_2; x_3, x_4) & = & \prod_{\sigma(3) < \sigma(4)} f(x_{\sigma(1)} , x_{\sigma
(2)}, x_{\sigma (3)} , x_{\sigma (4)})^{- \epsilon (\sigma)} \nonumber\\
\end{eqnarray}
There is a unique condition under which the two pairs of partial group
laws defining these two weak biextensions are commutative, thereby ensuring that both
$(g,h)$ and $(\gamma , \eta)$ define genuine biextensions. This condition is given by the vanishing of
the alternating map 
\[
l(x_1,x_2,x_3,x_4) = \prod_{\sigma \in \Sigma_4} f(x_{\sigma(1)} , x_{\sigma (2)} , x_{\sigma (3)} ,
x_{\sigma (4)})^{\epsilon ( \sigma)}
\]
determined by evaluating the 4-cocycle $f$ on the Pontrjagin product  $x_1 \cdot x_2 \cdot x_3
\cdot x_4$ of the four classes $x_i \in H_1(B) = B$. By proposition \ref{pr:altern}, both biextensions
$(g,h)$ and
$(\gamma ,
\eta)$ are then  alternating.  As explained in the
discussion following corollary
\ref{cor:alt}, they may therefore be described in cocyclic terms by the triples $(g,h,1)$ and
$(\gamma ,
\eta, 1)$ (in which the term $1$, which describes the alternating structure, is the trivial map
$1: B \la A $ sending every element of $B$ to the identity element of $A$). The compatibility
condition  (\ref{tri1}) between the first and second alternating structure is automatically satisfied
here. Since the  sections
$s^1_{x,y,z}$  and
$s^2_{x,y,z}$  satisfying the corresponding relations (\ref{phitriv}) are described by trivial maps,
this is also the case for  the  compatibility condition (\ref{t4}). 

\medskip

We summarize the previous discussion by the following 

\begin{proposition} Let $\cc$ be a monoidal category defined by a 4-cocycle $f(x_1 , x_2, x_3,
x_4)$, and for which the condition 
\[
\prod_{\sigma \in \Sigma_4} (f(x_{\sigma(1)} , x_{\sigma (2)}, x_{\sigma (3)} , x_{\sigma
 (4)}))^{\epsilon (\sigma)} = 1
\]
is satisfied. The pair of triples $(g(x_1, x_2; x_3 \mid_2 x_4),\:  h(x_1;  x_2, x_3 \mid_2 w) , \:1)$
and
 $\mbox{($\gamma(x_1 \nolinebreak \mid_2 \nolinebreak x_2, x_3 ; x_4),$} \linebreak \mbox{$\eta
(x_1 \mid_2 x_2; x_3, x_4),
\: 1$})$ defined by setting  
 
\begin{eqnarray} 
g(x_1, x_2; x_3 \mid_2 x_4) & = &\prod_{\sigma(1) < \sigma(2)}  f(x_{\sigma(1)} , x_{\sigma (2)},
x_{\sigma (3)} , x_{\sigma (4)})^{- \epsilon (\sigma)}\\
 h(x_1;  x_2, x_3 \mid_2 x_4) & = &  \prod_{\sigma(2) < \sigma(3)} f(x_{\sigma(1)} , x_{\sigma (2)},
x_{\sigma (3)} , x_{\sigma (4)})^{\epsilon (\sigma)}  \nonumber\\
\gamma(x_1 \mid_2 x_2, x_3 ; x_4) & = & \; \; \; \; h(x_1;  x_2, x_3 \mid_2 x_4)\nonumber \\
\eta (x_1\mid_2 x_2; x_3, x_4) & = & \prod_{\sigma(3) < \sigma(4)} f(x_{\sigma(1)} , x_{\sigma (2)},
x_{\sigma (3)} , x_{\sigma (4)})^{- \epsilon (\sigma)} \nonumber
\end{eqnarray}
respectively determine, for every fixed $x_4 \in B $ and 
 every fixed
$x_1 \in B$, a biextension structure. Together they  define an alternating triextension of $B
\times B \times B$ by
$A$. 
\end{proposition}

Suppose now that  this alternating triextension is trivial.  The trivializing section $\sigma_{x,y,z}$
of
$E_{x,y,z}$  (definition \ref{def:triex}) is then described by a map
$\theta: B^3 \la A$  such
that the  equations
\begin{eqnarray*}
g(x_1, x_2; x_3 \mid_2 x_4) & = & \frac{\theta(x_1 + x_2, x_3, x_4)}{\theta(x_1 , x_3, x_4)\:
\theta(x_2 , x_3, x_4)} \\
 h(x_1;  x_2, x_3 \mid_2 x_4) & = & \frac{\theta(x_1, x_2 + x_3, x_4)}{\theta(x_1, x_2, x_4)\:
\theta(x_1, x_3, x_4)} \\
\eta (x_1\mid_2 x_2; x_3, x_4) & = & 
\frac{\theta(x_1, x_2, x_3 + x_4)}{\theta (x_1, x_2, x_3) \: \theta (x_1, x_2, x_4)}\\
\theta (x,x,z) & = & 1\\
\theta (x,z,z) & = & 1
\end{eqnarray*}
are satisfied.  In that case the weak (2,2)-extension (\ref{2mock}) determined by 
proposition \ref{weak22ext} is a genuine (2,2)-extension. We may now introduce  additional cochains
$\theta (x_1, x_2 \mid_2 x_3)$ and $\theta (x_1 \mid_2 x_2 , x_3)$   defined by
\begin{eqnarray*}
\theta (x_1, x_2 \mid_2 x_3) & = & \theta (x_1, x_2, x_3) \\
\theta (x_1 \mid_2  x_2, x_3) & = & \theta (x_1, x_2, x_3)
\end{eqnarray*}
These  respectively describe the commutativity isomorphisms in the strict Picard categories
$\mathcal{E}_{ (\: , x_4)}$ and  $\mathcal{E}_{(x_1, \:)}$, so that the pair $(\phi (x_1,x_2, x_3
\mid_2 x_4), \theta (x_1, x_2 \mid_2 x_4))$ for
$x_4$ fixed and the pair $(\phi (x_1 \mid_2 x_2, x_3, x_4), \theta  (x_1 \mid_2 x_2, x_4))$ for
$x_1$ fixed each satisfy the  cocycle conditions (\ref{eq:braid})-(\ref{eq:ps}).  This
(2,2)-extension is automatically alternating, as may be verified by a discussion parallel to that of 
\S
\ref{par: altcom}, or by a cocyclic argument.

\bigskip

 While we could pursue this analysis in cocyclic
terms of the (2,2)-extension  $\mathcal{E}$, it is more
expedient to return to a 2-categorical framework. A trivialization of 
 $\mathcal{E}$ consists of a trivialization,
for each $x$ ({\it resp.} each
$w$) in $B$ of the Picard stack
$\mathcal{E}_{(x, \:)}$ ({\it resp.}
$\mathcal{E}_{(\: , w)}$), together with a compatibility condition between these trivializations.
Returning to the definition (\ref{2mock}) of
$\mathcal{E}$, we see that such a trivialization of $\mathcal{E}_{(x, \:)}$ (compatible with the
Picard structure) consists, once a family of choices of 1-arrows 
\be
\label{rxy}
R_{x,y}: X_y X_x  \la  X_x X_y
\ee 
have been made\footnote{after a preliminary choice of a  representative object $X_x$ of the
fixed object
$x$, and a family of representatives $X_y$ of the varying $y \in Y$.},  in a ``hexagon"  2-arrow
\[
H_{[x \mid_2 y,z]}: R_{x,y} \circ R_{x,y'} \Longrightarrow R_{x,yy'}
\]
 in $\cc$.  Compatibility of this trivialization with the associativity isomorphism in
$\mathcal{E}_{(x,
\:)}$ implies that this hexagon 2-arrow satisfies the axiom denoted $(\bullet \otimes (\bullet \otimes
\bullet
\otimes \bullet))$  in \cite{kv:2-cat}, and which we already encountered in a somewhat
different context. A trivialization of $\mathcal{E}_{(\: , w)}$ similarly
defines the ``hexagon" 2-arrow between the 1-map defined as in 
(\ref{u+1v}) (for $Y= X_w$) from the 1-arrow obtained by composing $R_{x,w}$ and $R_{x',w}$ and  the
1-arrow $R_{xx',w}$, and this then satisfies the corresponding axiom  $((\bullet \otimes \bullet
\otimes \bullet)
\otimes \bullet))$. The compatibilities of these 2-arrows with the commutativity
isomorphisms determined by the strict Picard structure yields the axioms $((\bullet \otimes \bullet)
\otimes (\bullet \otimes \bullet))$  on the 2-category
$\cc$. Finally, the  compatibility of this pair of hexagon 2-arrows with each other implies
that the 2-category $\cc$ is endowed with a slight variant of Kapranov-Voevodsky's 2-braiding
axioms, which we called a Z-braiding in \cite{lb:2-gerbes} chapter 8   (see also 
\cite{ba-ne:higher} for a discussion of this supplementary axiom). 

\medskip
If we now require that the chosen trivalization (\ref{rxy}) of the
(2,2)-extension $\mathcal{E}$ to be compatible with its anti-symmetry structure, defined as in
lemma \ref{lem:antisym}, we must further require that  there exists for all  $x,y \in B$ a 2-arrow 

\be
\label{ar:sxy}
\xymatrix{
{1_{X_yX_x}} \ar@{=>}[r]^>>>>>{S_{x,y}} & {R_{y,x} \circ R_{x,y}}
}
\ee
with source the identity 1-arrow
 on $X_yX_x$. This 2-arrow
automatically satisfies the two conditions which define on $\cc$ the structure of
 a strongly braided 2-category \cite{lb:2-gerbes} (in other words what J. Baez 
calls a strongly involutory monoidal category \cite{ba:2Hilb}).  Finally the compatibility of the
trivialization of
$\mathcal{E}$ with its alternating (rather than simply
 anti-symmetric) structure  manifests itself in a trivialization of the Picard
stack
$\Delta
\mathcal{E}$ obtained, as required in definition
\ref{22alt}, by restricting
$\mathcal{E}$ above the diagonal. Since the alternating structure on $\mathcal{E}$ is determined by
the identity arrow $1_{X_xX_x} \in (\Delta \mathcal{E})_{x}$, this compatibility may be interpreted as
a 2-arrow
\[
\xymatrix{
{1_{X_xX_x}} \ar@{=>}[r]^{S_x} & {R_{x,x}}}
\]
 By compatibility of this 2-arrow with the group law on $\Delta \mathcal{E}$, $S_x$ is
 additive in $x$. Furthermore, the required compatibility of its square with the
trivialization of $\mathcal{E}^2$ determined by the anti-symmetry structure on $\mathcal{E}$ is the
assertion that the composite 2-arrow in diagram (8.4.8) of 
\cite{lb:2-gerbes} coincides with our 2-arrow $S_{x,x}\:$ (\ref{ar:sxy}). Observe also that the
compatibility condition mentioned  in \cite{lb:2-gerbes} (8.4.6) is in fact a consequence of the
required  additivity in
$x$ of
$S_x$, and therefore must not be imposed here as a supplementary condition. A trivialization of
the alternating $(2,2)$-extension $\mathcal{E}$  thus determines  on
$\cc$  what we have called a strictly symmetric monoidal  2-category structure. These  strictly
symmetric structures on monoidal  2-categories with associated groups $B$ and $A$ are classified 
by the group
$\mathrm{Ext}^3(B, A)$. In the 2-stack case, this is a genuine invariant, which classifies these
structures up to equivalence. On the other hand, in the 2-category case, the vanishing of this
group is automatic, since it is a higher  $\mathrm{Ext}$ group in the category of
abelian groups. The 2-category
$\cc$ is  therefore equivalent to the trivial one with invariant groups $B$ and
$A$.  In particular, forgetting all the symmetry stucture, this  implies that  the underlying monoidal
2-category
$\cc$ is equivalent to the trivial one, and  the original
4-cocycle
$f(x,y,z,w)$ which defined $\cc$ is  then cohomologous to zero.

\end{document}